\documentclass{msml2020} 

\usepackage{bm,amsmath}
\usepackage{mathtools,empheq,algorithm,algorithmic}
\usepackage{etoolbox,booktabs}
\newcommand{\algorithmicdoinparallel}{\textbf{do in parallel}}
\makeatletter
\AtBeginEnvironment{algorithmic}{%
  \newcommand{\FORALLP}[2][default]{\ALC@it\algorithmicforall\ #2\ %
    \algorithmicdoinparallel\ALC@com{#1}\begin{ALC@for}}%
}
\makeatother
\mathtoolsset{showonlyrefs}

\DeclareMathOperator*{\argmin}{arg\,min}

\newcommand{\cN}{\mathcal{N}}

\newcommand{\bX}{\bm{X}}
\newcommand{\bx}{\bm{x}}

\newcommand{\bZ}{\bm{Z}}
\newcommand{\bz}{\bm{z}}
\newcommand{\bW}{\bm{W}}
\newcommand{\balpha}{\bm{\alpha}}

\newcommand{\mc}[1]{\mathcal{#1}}

\newcommand{\EE}{\mathbb{E}}

\newcommand{\RR}{\mathbb{R}}

\newcommand{\inv}{^{\raisebox{.2ex}{$\scriptscriptstyle-1$}}}
\newcommand{\oon}{\frac{1}{N}}

\newcommand{\dij}[1]{\delta_{#1}}
\newcommand{\place}[2]{%
	\underset{\substack{\uparrow\\#1}}{#2}%
}

\newcommand{\transpose}{^{\operatorname{T}}}

\newcommand{\figwidth}{{0.95\textwidth}}

\newcommand{\half}{\frac{1}{2}}

\newcommand{\eps}{\epsilon}

\newcommand{\abs}[1]{\left|#1\right|}

\newcommand{\ud}{\,\mathrm{d}}

\newtheorem{theo}{Theorem}[section]

\newtheorem{rem}[theo]{Remark}

\newtheorem{defi}[theo]{Definition}

\title[Deep Fictitious Play for Finding Markovian Nash Equilibrium]{Deep Fictitious Play for Finding Markovian Nash Equilibrium in Multi-Agent Games}

\usepackage{times}

 \msmlauthor{\Name{Jiequn Han} \Email{jiequnh@princeton.com}\\
 \addr Department of Mathematics, Princeton University, Princeton, NJ 08544-1000, USA
 \AND
 \Name{Ruimeng Hu} \Email{rh2937@columbia.edu}\\
 \addr Department of Statistics, Columbia University, New York, NY 10027-4690, USA
 }

\begin{document}
\maketitle

\begin{abstract}
We propose a deep neural network-based algorithm to identify the Markovian Nash equilibrium of general large $N$-player stochastic differential games. Following the idea of fictitious play, we recast the $N$-player game into $N$ decoupled decision problems (one for each player) and solve them iteratively. The individual decision problem is characterized by a semilinear Hamilton-Jacobi-Bellman equation, to solve which we employ the recently developed deep BSDE method. The resulted algorithm can solve large $N$-player games for which conventional numerical methods would suffer from the curse of dimensionality. Multiple numerical examples involving identical or heterogeneous agents, with risk-neutral or risk-sensitive objectives, are tested to validate the accuracy of the proposed algorithm in large group games. Even for a fifty-player game with the presence of common noise, the proposed algorithm still finds the approximate Nash equilibrium accurately, which, to our best knowledge, is difficult to achieve by other numerical algorithms. 
\end{abstract}

\begin{keywords}
stochastic differential game, Markovian Nash equilibrium, fictitious play, deep learning
\end{keywords}

\section{Introduction}


{The development of communication technologies makes the world more interactive than ever, and crucial decisions are generally made in a common environment nowadays where multiple players/agents\footnote{Hereafter, we shall use \textit{player} and \textit{agent} interchangeably.} act reciprocally. The stochastic differential game provides an elegant framework for studying such scenarios and has been widely used in the areas of financial mathematics, economics, management science, engineering, etc. For instance, under this framework, \cite{CaFoSu:15} analyze systemic risk caused by inter-bank lending and borrowing,  \cite{HACT2017} study the income and wealth distribution in macroeconomics,
and \cite{BeSiYaYa:14} model the reinsurance games between insurance companies, just to list a few. }


Among all applications of stochastic differential games, a core problem is to find the Nash equilibrium. A Nash equilibrium refers to a set of all players' strategies in which nobody has the incentive to deviate, given the others' strategies fixed, and a Markovian strategy means the decision is made based only on the current situation of all players. 
Under the Markovian setting in continuous time, the Nash equilibrium is characterized by $N$-coupled Hamilton-Jacobi-Bellman (HJB) equations, each of which is usually high-dimensional when $N$ is large.
In some special settings, one can find the closed-form Nash equilibrium. However, an analytic solution is generally unavailable, and one has to resort to numerical solutions. 
The notorious computational difficulty ``curse of dimensionality" is  often encountered in this situation due to the high dimensionality of state variables.
So it is crucial but challenging to develop numerical algorithms to accurately and efficiently find the Nash equilibrium in large $N$-player stochastic differential games.


In recent years the deep neural network (DNN) has shown its remarkable ability in representing and approximating high-dimensional functions in numerical computation, ranging from artificial intelligence (cf. \cite{Be:09,LeBeHi:15}) such as computer vision and speech recognition, to scientific computing (see e.g., \cite{CaTr:17,GaDu:17,HaJeE:18,ZhHaE:18,han2019uniformly}).
Considering the astonishing performance of the DNN, we shall leverage it to help identify the Nash equilibrium through solving high-dimensional coupled HJB equations. To achieve this, we design a scalable algorithm, which we refer to as \emph{deep fictitious play}, based on a combination of DNN and the idea of fictitious play.
Fictitious play is an iterative scheme in which, at each stage, the whole coupled game is decomposed to independent decision problems for each player by considering her opponents' strategies are fixed and following their past play. Such a scheme enables the algorithm to be decoupled and parallel between players. The individual decision problem is then translated into a single HJB equation and solved via the deep BSDE method (\cite{EHaJe:17,HaJeE:18}). Since decision problems are solved individually for all players through general HJB equations, the algorithm can accommodate games with arbitrary heterogeneity, and with risk-neutral or risk-sensitive cost functional. 
Note that the idea of deep fictitious play is recently proposed in \cite{Hu2:19} as well for solving open-loop Nash equilibrium. However, due to the entirely different strategy set, the algorithm developed in \cite{Hu2:19} is inapplicable for Markovian Nash equilibrium (see Remark~\ref{rem_lastvsaverage} for details).

An alternative approach to studying large $N$-player game is mean-field game theory, introduced by \cite{LaLi1:2006,LaLi2:2006,LaLi:2007} and by \cite{HuCaMa:07,HuMaCa:06} around the same time. In mean-field game theory, the limiting solution as $N \to \infty$ is used to approximate the game of finite population. Under mild assumptions, \cite{CaDe:13} have shown that the approximation error is of order $N^{-1/(d+4)}$, where $d$ is the dimension of each player's state variable. The limitation of the mean-field game is threefold. First, the approximation error is not well controlled for a moderate $N$. Second, it requires lots of symmetry. To make the mean-field approximation work, the players need to be indistinguishable, i.e., they control their states in the same way and have identical objectives. Third, the numerical computation of the mean-filed game is usually infeasible when there is common noise in addition to idiosyncratic noise. In contrast to these limitations, as shown in the description of the algorithm and numerical results, the proposed algorithm can solve general $N$-player stochastic differential games with the presence of heterogeneity and common noise.

Another line of research relevant to this work is the so-called multi-agent reinforcement learning (MARL) (see~\cite{Busoniu08} for an earlier review), which recently has gained increasing research attention (see e.g.~\cite{foerster2016learning,zhang2018fully}), given the remarkable success of reinforcement learning in many applications, especially accompanied with the powerful function approximation based on DNN. The main difference is that MARL is usually model-free and in discrete time while our work considers model-based games in continuous time, for which we can leverage the HJB equation and BSDE reformulation to develop efficient algorithms. 

In short, our algorithm focuses on the moderate number of the agents, i.e., $10 \leq N \leq 100$, in which case conventional numerical methods lose their efficiencies while the mean-field theory has not yet been accurate in approximating the finite population game. The strength of the proposed algorithm is the ability to deal with heterogeneous agents, to which scenarios mean-field games are not good at in general, and the ease of dealing with common noise.
The rest of the paper is organized as follows. Section~\ref{sec_math} gives the mathematical formulation of general $N$-player stochastic differential games in continuous time. The algorithm based on fictitious play and DNN is detailed in Section~\ref{sec_method}, followed by numerical examples in Section~\ref{sec_numerics}. We conclude in Section~\ref{sec_conclusion} and provide technical details in appendices.

\section{Mathematical Formulation}\label{sec_math}
We consider $N$-player non-cooperative stochastic differential games described by the following dynamics,
\begin{equation}\label{def_Xt}
\ud \bm X_t^{\bm \alpha} = b(t, \bm X_t^{\bm \alpha}, \bm \alpha(t, \bm X_t^{\bm \alpha})) \ud t + \Sigma(t, \bm X_t^{\bm \alpha}) \ud \bm W_t, \quad \bm X_0 = \bm x_0,
\end{equation}
where $\bm X_t^{\bm \alpha}$
is an $\RR^n$-valued \emph{common} state process influenced by Markovian controls $\bm \alpha = (\alpha^1, \ldots, \alpha^N)$, a collection of all players' strategies. Each $\alpha^i$, as the control of player $i$, is a Borel measurable function $[0,T] \times \RR^n \hookrightarrow \mc{A}^{i} \subset \RR^{n_i}$. 
$b$ and $\Sigma$ are deterministic functions denoting the drift and diffusion coefficients of the common state, $b\colon [0,T] \times \RR^n \times \mc{A}  \hookrightarrow \RR^n$,  $\Sigma\colon [0,T] \times \RR^n \hookrightarrow \RR^{n \times k}$, where $\mathcal{A} =  \otimes_{i=1}^N \mc{A}^i$ is the space for joint control $\balpha$, and $\bm W$ is a $k$-dimensional standard Brownian motion on a filtered probability space $(\Omega,\mathbb{F},\{\mathcal{F}_t\}_{0 \le t \le T}, \mathbb{P})$.

Denote by $\mathbb{A}^i$ the set of admissible $\mc{A}^i$-valued strategy for player $i \in \mc{I} := \{1, 2, \ldots, N\}$, and by $\mathbb{A} = \otimes_{i=1}^N \mathbb{A}^i$ the product space of $\mathbb{A}^i$. Given the other's strategy $\alpha^{j} \in \mathbb{A}^{j}$, $j \neq i$,  player $i$ chooses $\alpha^i$ to minimize the cumulative cost $J_0^i$, 
\begin{equation}\label{def_J}
J_t^i(\bm \alpha) := \EE\left[\int_t^T f^i(s, \bm X_s^{\bm \alpha}, \bm \alpha(s,\bm  X_s^{\bm \alpha})) \ud s + g^i(\bm X_T^{\bm \alpha}) \right],
\end{equation}
or its risk-sensitive version
\begin{equation}\label{def_J'}
J^i_t(\bm \alpha) := \EE\left[\theta_i \exp\left\{\theta_i\left(\int_t^T f^i(s, \bm X_s^{\bm \alpha}, \bm \alpha(s, \bm X_s^{\bm \alpha}) \ud t + g^i(\bm X_T^{\bm \alpha})\right)\right\} \right],
\end{equation}
where the running cost $f^i\colon [0,T] \times \RR^{n} \times \mathcal{A} \hookrightarrow \RR$ and the terminal cost $g^i\colon \RR^{n} \hookrightarrow \RR$ are deterministic measurable functions.

We are interested in solving the above game in terms of finding the Nash equilibrium, in particular, the  Markovian Nash equilibrium. 

\begin{defi}
A Markovian Nash equilibrium is a tuple $\bm \alpha^\ast = (\alpha^{1, \ast}, \ldots, \alpha^{N, \ast}) \in \mathbb{A}$ such that 
\begin{equation}
\forall i \in \mc{I}, \emph{ and } \alpha^i \in \mathbb{A}^i, \quad J^i_0(\bm \alpha^\ast) \leq J^i_0(\alpha^{1, \ast}, \ldots, \alpha^{i-1, \ast}, \alpha^i, \alpha^{i+1, \ast}, \ldots, \alpha^{N, \ast}).
\end{equation}
\end{defi}

The above problem \eqref{def_Xt}--\eqref{def_J'} is more general than the usual finite-player stochastic differential games or mean-field games in two aspects. Firstly, in addition to the usual (risk-neutral) cost functional \eqref{def_J}, we also consider a risk-sensitive version \eqref{def_J'}, where $\theta_i$ is a parameter characterizing how risk-averse/seeking player $i$ is. This 
flexibility allows us to model much broader classes of games that accommodate the players' attitudes to risk.
Secondly and more importantly, $\bm X_t$ in \eqref{def_Xt} is a common state process that is influenced by all players, as opposed to the traditional case that player $i$ can only control her private state. The former feature is common in the literature of economics (see e.g., \cite{DoJoVaSo:00,PrSe:04,Va:11}), thus, we think it is important to include it in \eqref{def_Xt}. However, this leads to a stronger coupling problem, which is unsurprisingly harder to deal with both theoretically and numerically. The difficulty even persists in the limiting problem as $N \to \infty$ with indistinguishable players, when allowing $\alpha^i$ entering into others' states. This is called the extended mean-field game and it has attracted lots of attention recently (see e.g., \cite{GoVo:13,GoPaVo:14,GoVo:16,CaLe:18}).

Note that by choosing $b$ and $\Sigma$ in \eqref{def_Xt} properly, one can reduce the formulation \eqref{def_Xt} to the simpler case where each player controls her \emph{private} state through $\alpha^i$. We highlight this relation in the following remark. 
\begin{rem}
Let $n = dN$ and $b^\ell \equiv b^\ell(t, \bm x, \alpha^i)$ for $\ell = (i-1)d+1, \ldots, id$,
then the problem \eqref{def_Xt}--\eqref{def_J'} is the usual modeling in financial mathematics literature, where the $i^{th}$ player's $d$-dimensional private state $(X_t^{(i-1)d+1}, \ldots, X_t^{id})$ is controlled by $\alpha^i$ only. A benchmark example we shall study in Section~\ref{sec_example1} is a model of inter-bank borrowing and lending proposed in \cite{CaFoSu:15}, where $d = 1$, and the dynamics for each player is
\begin{equation}
\ud X_t^i = [a(\overline X_t - X_t^i)   + \alpha_t^i ] \ud t  + \sigma \left(\rho \ud W_t^0 + \sqrt{1-\rho^2} W_t^i\right), \quad \overline X_t = \frac{1}{N}\sum_{i=1}^N X_t^i.
\end{equation}

\end{rem}

In the Markovian setting, finding a Nash equilibrium is related to solving $N$-coupled HJB equations. To this end, we define the value function of player $i$ by
\begin{equation}
V^i(t, \bm x) = \inf_{\alpha^i \in \mathbb{A}^i} \EE[J_t^i(\bm \alpha) \vert \bm X_t = \bm x].
\end{equation}
Using the dynamic programming principle, the HJB system reads
\begin{align}
\label{def_HJB}
\begin{dcases}
V_t^i + \inf_{\alpha^i \in \mc{A}^i} G^i(t,\bm x,\bm \alpha, \nabla_{\bx} V^i, V^i)  + \half \text{Tr}(\Sigma\transpose \text{Hess}_{\bx} V^i \Sigma) = 0,\\
V^i(T,\bx) = g^i(\bx), \quad i \in \mc{I},
\end{dcases}
\end{align}
where $G^i$ is given by 
\begin{align}
&\text{under objective \eqref{def_J}}, \quad G^i = G^i(t,\bm x, \bm \alpha, \bm p) =  b(t,\bm x, \bm \alpha) \cdot \bm p + f^i(t, \bm x, \bm \alpha), \\
&\text {under  objective \eqref{def_J'}}, \quad  G^i = G^i(t,\bm x, \bm \alpha, \bm p,s) =  b(t,\bm x, \bm \alpha) \cdot \bm p + \theta_i sf^i(t, \bm x, \bm \alpha).  
\end{align}
Here $\nabla_{\bx} V,~\text{Hess}_{\bx} V$ denote 
the gradient and the Hessian of function $V$ respect to $\bm x$ and $\text{Tr}$ denotes the trace of a matrix. In the risk-neutral case \eqref{def_J}, $G^i$ is in fact the Hamiltonian and is usually denoted by $H^i$. Nevertheless, to unify the notation between two objectives \eqref{def_J} and \eqref{def_J'}, we shall stick to $G^i(t, \bm x, \bm \alpha, \bm p, s)$ for the rest of the paper, although it does not depend on $s$ under objective \eqref{def_J}. Note that while minimizing $\alpha^i$ in the equation for $V^i$ in \eqref{def_HJB}, the policies $(\alpha^1, \dots,\alpha^{i-1}, \alpha^{i+1}, \dots, \alpha^N)$ are given and fixed. In other words, $V^i$ implicitly depends on the other players' strategies, thus on $V^j$.

Throughout the paper, we assume that there is a unique classical solution to the HJB system \eqref{def_HJB}. Moreover, we require that the minimizer $\argmin_{\alpha^i \in \mc{A}^i} G^i(t, \bm x, \bm \alpha, \bm p, s)$ exists, is unique and explicit in other arguments, $\forall$  $i \in \mathcal{I}$, $(t, \bm x, \bm p, s)$ and $\alpha^j \in \mc{A}^j$ with $j \neq i$.

\section{Methodology}
\label{sec_method}
\subsection{Fictitious Play}
Fictitious play was firstly introduced by Brown in the static game \cite{Br:49,Br:51}, and was recently adapted to the mean-field setting by \cite{CaHa:17,BrCa:18}. It is a simple yet important learning idea in game theory. The key is to decouple the $N$-player game into $N$ individual decision problems where opponents' strategies are fixed and assumed to follow their past play. These $N$ individual problems are solved iteratively, starting from stage 1. At stage $m$, we assume that the opponents' strategies are their stage $(m-1)$'s best responses. To better describe this procedure mathematically, we first summarize the notations that shall be used below.

\begin{itemize}
	\item $\mc{A}^i \subset \RR^{n_i}$, the range of player $i$'s strategy $\alpha^i$.  $\mc{A} = \otimes_{i=1}^N \mc{A}^i$, the control space for all players, and $\mc{A}^{-i} = \otimes_{j\neq i} \mc{A}^j$, the control space for all players but $i$. The same applies to the admissible spaces of measurable functions $\mathbb{A}^i$, $\mathbb{A}$ and $\mathbb{A}^{-i}$. Note that, with these notations, $\inf_{\alpha^i\in \mc{A}^i}$ and $\inf_{\alpha^i\in \mathbb{A}^i}$ means seeking for an optimal vector in $\mc{A}^i$ and $\mc{A}^i$-valued function (strategy), respectively.  
	\item $\bm{\alpha} = (\alpha^1, \alpha^2, \ldots, \alpha^N)$, a collection of all players' strategy profiles. With a negative superscript, $\bm{\alpha}^{-i} = (\alpha^1, \ldots, \alpha^{i-1}, \alpha^{i+1}, \ldots, \alpha^{N})$ means  the strategy profiles excluding player $i$'s. If a non-negative superscript $m$ appears, $\bm{\alpha}^m$ is a $N$-tuple standing for the strategies of all players at stage $m$. When both exist,  $\bm{\alpha}^{-i,m} =  (\alpha^{1,m}, \ldots, \alpha^{i-1,m}, \alpha^{i+1,m}, \ldots, \alpha^{N,m})$ is a $(N-1)$-tuple representing strategies excluding player $i$ at stage $m$. 
\end{itemize}

Assume that we start with a guess of the solution $\bm\alpha^0 \in \mathbb{A}$. The idea of fictitious play motivates us to consider an iterative algorithm according to the following rules.
At stage $m+1$, $\bm \alpha^m$ is observed by all players, and player $i$'s decoupled decision problem is 
\begin{equation}\label{def_J_fictitious}
\inf_{\alpha^i \in \mathbb{A}^i} J_0^i(\alpha^i; \bm \alpha^{-i,m}),
\end{equation}
where $J^i_0$ is defined in \eqref{def_J} or \eqref{def_J'}, and the state process $\bm X_t$ is given in \eqref{def_Xt} with $\bm \alpha$ being replaced by $(\alpha^i, \bm \alpha^{-i,m})$. The optimal strategy, if ever exists, is denoted by $\alpha^{i, m+1}$.
The problem \eqref{def_J_fictitious} for all $i \in \mc{I}$ are solved simultaneously using $\bm \alpha^{-i, m}$, and the optimal responses together form $\bm \alpha^{m+1}$. Due to the Markovian structure, the problem \eqref{def_J_fictitious} is translated into a HJB equation
\begin{equation}
\label{eq:FP_PDE}
V_t^{i,m+1} + \inf_{\alpha^i \in \mc{A}^i} G^i(t, \bm x, (\alpha^i, \bm \alpha^{-i,m}(t, \bm x)), \nabla_{\bx} V^{i, m+1}, V^{i, m+1}) + \half \text{Tr}(\Sigma\transpose \text{Hess}_{\bx} V^{i, m+1} \Sigma) = 0,
\end{equation}
with the terminal condition $V^{i,m+1}(T,\bx) = g^i(\bx)$.
If the classical solution ever exists, the optimal strategy at stage $m+1$ for player $i$ is given by
\begin{equation}\label{def_alphaast}
\alpha^{i, m+1}(t, \bm x) = \argmin_{\alpha^i \in \mc{A}^i} G^i(t, \bm x, (\alpha^i, \bm \alpha^{-i,m}(t, \bm x)), \nabla_{\bx} V^{i, m+1}(t, \bx), V^{i, m+1}(t, \bx)).
\end{equation}
Solving~\eqref{eq:FP_PDE} for all $i\in \mc{I}$ completes one stage in the loop of fictitious play.

Given the iterative procedure described above, one can naturally ask: (1) does $\alpha^{i,m}$ always exist; (2) if yes, does $\bm \alpha^m$ converge; and (3) if yes, is the limiting strategy $\bm \alpha^\infty$ admissible and does it form a Nash equilibrium.
The first question closely depends on the choice of the initial belief $\bm \alpha^0$. Usually a regular enough $\bm \alpha^0$ (plus regular $b, \Sigma, f^i, g^i$) will ensure the unique classical solution $V^{i,1}$, thus ensure the existence of $\bm \alpha^1$. Then it is likely that the regularity of $\bm \alpha$ persists from stage to stage. 
The remaining two questions are tough in general.
For the second question, there is no universal criterion that can guarantee the convergence.
Even in the static finite-action games, there are numerous examples when it converges (\cite{Ro:51, Mi:61, MiRo:91,MoSh:96,MoSh2:96, Be:05, HoSa:02}) and when it does not (\cite{Sh:64, Jo:93,MoSe:96,FoYo:98,KrSj:98}). The third question is related to the stability of the game, and has to be analyzed case by case.

Instead of answering the theoretical questions raised above, this paper focuses on a practical numerical scheme for finding the Markovian Nash equilibrium, especially when $N$ is large. As proof of methodology, we shall present four examples to show the performance of the proposed algorithm. These examples cover a large variety of stochastic differential games, including homogeneous/heterogeneous, risk-neutral/risk-sensitive ones.  We remark that the first three examples are chosen so that they have closed-form solutions or can be solved via other numerical approaches in low dimensions, such that we can benchmark numerical solutions. The algorithm in fact can be applied in much more general games, as demonstrated in the last example.

\begin{rem}
The learning procedure described above is slightly different from the original version proposed by \cite{Br:49,Br:51} where two-player normal-form games are studied. Targeting at pure or mixed Nash equilibrium, the player therein chooses the best response against the empirical distribution of her opponent's past play. That is, if her opponent uses strategies $\alpha_1, \ldots, \alpha_k$ during the first $k$ stages, where $\alpha_i$ is from a finite strategy set, then at stage $k+1$, her response will be a pure strategy that maximizes her expected payoff with respect to her opponent's mixed strategy $\frac{1}{k}\sum_{i=1}^k \delta_{\alpha_i}$.

In this paper, we prefer to use the last stage information, instead of the average of the past. The reason will be explained after we introduce our deep learning algorithm, see Remark~\ref{rem_lastvsaverage}.

\end{rem}

\subsection{A Deep Learning Algorithm Based on Fictitious Play}\label{sec_algorithm}
In order to find the Markovian Nash equilibrium through fictitious play, the main challenge is numerically solving the PDE \eqref{eq:FP_PDE} at each stage. When the dimension $d$ is large, conventional numerical algorithms soon lose their efficiency. Here we employ the deep BSDE method in \cite{EHaJe:17,HaJeE:18} to deal with the high-dimensionality. 

Instead of solving the parabolic PDE \eqref{eq:FP_PDE}, the deep BSDE method recasts the problem into an optimization problem based on the associated BSDE. Under the standing assumptions, the function $\alpha^{i, m+1}(t, \bm x, \bm p, s; \bm \alpha^{-i,m}) = \argmin_{\alpha^i\in\mc{A}^i} G^i(t, \bm x, (\alpha^i, \bm \alpha^{-i,m}(t, \bm x)), \bm p, s)$ is unique and explicit. We plug it
into \eqref{eq:FP_PDE} and rewrite the equation as a semilinear parabolic PDE:
\begin{multline}
     V_t^{i,m+1} + \half \text{Tr}(\Sigma\transpose \text{Hess}_{\bx} V^{i,m+1} \Sigma) + \mu^i(t, \bx; \balpha^{-i, m})\cdot \nabla_{\bx} V^{i, m+1} \\ +h^i(t, \bx,  V^{i, m+1}, \Sigma\transpose\nabla_{\bx} V^{i, m+1}; \balpha^{-i, m})=0.
\end{multline}
Note that we replace $\nabla_{\bx} V^{i, m+1}$ by $\Sigma\transpose\nabla_{\bx} V^{i, m+1}$ in the term $h^i$ for the sake of simplicity when describing the algorithm, and that we treat $\balpha^{-i,m}$ as known functions of the PDE due to its exogeneity. The concrete equations of the considered numerical examples are presented in the associated subsections or Appendix~\ref{appendix_PDE}.
For less cumbersome notations, we drop the superscript $m$ that denotes the index of fictitious play and consider the PDE
\begin{align}
\label{eq:FP_PDE_explicit}
\begin{dcases}
V_t^i + \half \text{Tr}(\Sigma\transpose \text{Hess}_{\bx} V^i \Sigma) + \mu^i(t, \bx; \balpha^{-i})\cdot \nabla_{\bx} V^{i}+h^i(t, \bx,  V^{i}, \Sigma\transpose\nabla_{\bx} V^{i};  \balpha^{-i})=0,\\
V^i(T,\bx) = g^i(\bx).
\end{dcases}
\end{align}
This PDE is intimately related to the following BSDE: 
\begin{empheq}[left=\empheqlbrace]{align}
    &\bX_t^i = \bx_0 + \int_{0}^{t}\mu^i(s,\bX_s^i; \balpha^{-i}(s, \bX_s^i))\, \mathrm{d}s + \int_{0}^{t}\Sigma(s,\bX_s^i)\, \mathrm{d}\bW_s, \label{eq:BSDE_forward} \\
    &Y_t^i = g^i(\bX_T^i) + \int_{t}^{T}h^i(s,\bX_s^i,Y_s^i,\bZ_s^i;  \balpha^{-i}(s, \bX_s^i))\, \mathrm{d}s - \int_{t}^{T}(\bZ_s^i)\transpose\, \mathrm{d}\bW_s, \label{eq:BSDE_backward}
\end{empheq}
where $\bW_t$ is a $k$-dimensional Brownian motion and $\bx_0$ is a square-integrable random variable independent of $\bW_t$.
Specifically, the nonlinear Feynman-Kac formula (cf. \cite{PaPe:92,ElPeQu:97,PaTa:99}) states that, under certain regularity conditions,
\begin{equation}
\label{eq:FK relation}
    Y_t^i = V^i(t, \bX_t^i) \quad\text{and}\quad \bZ_t^i = \Sigma(t, \bX_t^i)\transpose\nabla_{\bx} V^i(t, \bX_t^i)
\end{equation}
defines the unique solution to the BSDE \eqref{eq:BSDE_forward}--\eqref{eq:BSDE_backward}. Accordingly, we consider the variational problem
\begin{align}
&\inf_{Y_0^i,\{\bZ_t^i\}_{0\le t \le T}} \EE|g^i(\bX_T^i) - Y_T^i|^2, \label{eq:variational_form}\\
&s.t.~~ \bX_t^i = \bx_0 + \int_{0}^{t}\mu^i(s,\bX_s^i; \balpha^{-i}(s, \bX_s^i))\, \mathrm{d}s + \int_{0}^{t}\Sigma(s,\bX_s^i)\, \mathrm{d}\bW_s, \notag \\
 & \qquad Y_t^i = Y_0^i - \int_{0}^{t}h^i(s,\bX_s^i,Y_s^i,\bZ_s^i; \balpha^{-i}(s, \bX_s^i))\, \ud s + \int_{0}^{t}(\bZ_s^i)\transpose\, \ud\bW_s,  \notag
\end{align}
where $Y_0^i$ is $\mathcal{F}_0$-measurable and square-integrable, and $\bZ_t^i$ is a $\mathcal{F}_t$-adapted square-integrable process. 
Eq. \eqref{eq:FK relation} gives a minimizer of the above problem since the loss function attains zero when it is evaluated according to \eqref{eq:FK relation}.
In addition, if the BSDE \eqref{eq:BSDE_forward}--\eqref{eq:BSDE_backward} is wellposed  (under some regularity conditions), the minimizer must exist and is unique.

The deep BSDE method builds on the temporal discretization of \eqref{eq:variational_form}. Given a partition $\pi$ of size $N_T$ on the time interval $[0, T], 0 = t_0 < t_1 < \ldots < t_{N_T} = T$, we consider the discretized version of \eqref{eq:variational_form} based on the Euler scheme (the subscript $t_k$ in $\bX, Y, \bZ$ has been replaced by $k$ for simplicity):
\begin{align}
&\inf_{\psi_0\in \cN_0^{i'},~\{\phi_k\in \cN_k^i\}_{k=0}^{N_T-1} } \EE|g^i(\bX_{N_T}^{i,\pi}) - Y_{N_T}^{i,\pi}|^2, \label{eq:disc_objective}\\
&s.t.~~ \bX_0^{i,\pi}=\bx_0, \quad Y_0^{i,\pi} =\psi_0(\bX_0^{i,\pi}), \quad \bZ_{k}^{i,\pi}=\phi_k(\bX_{k}^{i,\pi}), \quad k=0,\dots,N_T-1\notag \\
&\qquad \bX_{k+1}^{i,\pi} = \bX_{k}^{i,\pi} + \mu^i(t_k,\bX_{k}^{i,\pi}; \balpha^{-i}(t_k, \bX_k^{i,\pi}))\Delta t_k +\Sigma(t_k,\bX_{k}^{i,\pi})\Delta \bW_{k}, \label{eq:disc_X_path} \\
 &\qquad Y_{k+1}^{i,\pi} = Y_{k}^{i,\pi} - h^i(t_k,\bX_{k}^{i,\pi},Y_{k}^{i,\pi},\bZ_{k}^{i,\pi}; \balpha^{-i}(t_k, \bX_k^{i,\pi}))\Delta t_k + (\bZ_{k}^{i,\pi})\transpose\Delta \bW_k, \label{eq:disc_Y_path}
\end{align}
where
  $\Delta t_k = t_{k+1} - t_k, \quad \Delta \bW_k = \bW_{ t_{k+1}}-\bW_{t_k}.$
Here $\cN_0^{i'}$ and $\{\cN_k^i\}_{k=0}^{N_T-1}$ are hypothesis spaces of player $i$ related to deep neural networks. The goal of the optimization is to find optimal deterministic maps $\psi_0^{i,\ast}, \{\phi_k^{i,\ast}\}_{k=0}^{N_T-1}$ such that the loss function is small. Intuitively, we expect that  \eqref{eq:disc_objective} defines a benign optimization problem close to \eqref{eq:variational_form} and $\psi_0^{i,*}, \{\phi_k^{i,*}\}_{k=0}^{N_T-1}$ provide us good approximations to $V^i(0, \cdot), \{\nabla_{\bx} V^i(t_k, \cdot)\}_{k=0}^{N_T-1}$, the solution of the original PDE \eqref{eq:FP_PDE_explicit}. 

In practice, the expectation in \eqref{eq:disc_objective} is approximated by the standard Monte Carlo sampling of \eqref{eq:disc_X_path}--\eqref{eq:disc_Y_path}. Given sample paths of $\{\bX_{k}^{i,\pi}\}_{k=0}^{N_T}$ and the associated white noises $ \{\bW_{t_k}\}_{k=0}^{N_T}$, one notes that the loss function \eqref{eq:disc_objective} can be interpreted as the final output of a very deep network after stacking all the subnetworks $\psi_0^i, \{\phi_k^i\}_{k=0}^{N_T-1}$ in sequence according to \eqref{eq:disc_Y_path}. This means that we can use backpropagation to derive the gradient of the loss function with respect to all the parameters in the neural networks and use stochastic gradient descent (SGD) to optimize all the parameters. We refer the interested readers to \cite{EHaJe:17,HaJeE:18,HaLo:18} for more detailed description and theoretical justification of the deep BSDE method.

\begin{rem}\label{rem_lastvsaverage}
According to the deep BSDE method, the best responses at stage $m$, $\bm \alpha^m(t, \bm x)$, are defined through the outputs of neural networks. Due to the direct feedback nature, $\bm X_t$ changes as strategies vary from stage to stage, so do the function evaluations $\bm \alpha^{m}(t, \bm X_t)$. Therefore, if one were interested in using $\overline{\bm \alpha}^{-i,m} := \frac{1}{m} \sum_{\ell =1}^ m \bm \alpha^{-i, \ell}$ to replace $\bm \alpha^{-i,m}$ in \eqref{def_J_fictitious}, the parameters of all neural networks from stage 1 to $m$ need to be saved, which is infeasible for the sake of computational memory.  This is indeed the very reason that we use $\bm \alpha^{-i, m}$ in objective \eqref{def_J_fictitious}.

Note that this is not the case for searching open-loop Nash equilibrium. There, $\bm \alpha^{m}$ are adapted processes to the filtration generated by Brownian motions, i.e., one can intuitively think $\bm \alpha^{m}$ as functions of $\bm W_{[0,T]}$, thus do not change as strategies/states vary, as long as the training paths of $\bm W$ are fixed. Therefore, using $\overline{\bm \alpha}^{-i,m}$ only requires constant memory, since one can update it by a weighted sum of $\overline{\bm \alpha}^{-i,m-1}$ and $\bm \alpha^{-i, m}$.

For further discussion on this discrepancy, we refer to a closely related work by \cite{Hu2:19} where a deep learning scheme is designed for open-loop Nash equilibrium based on fictitious play. Indeed, it is empirically observed in \cite{Hu2:19} that using $\overline{\bm \alpha}^{-i,m}$ leads to faster convergence in the open-loop case. We would expect the similar phenomenon in the Markovian case, if we were able to record all the best responses $\bm{\alpha}^m(t, \bm x)$.
However, as discussed above, due to the direct feedback nature, merely adapting the algorithm in \cite{Hu2:19} will not be efficient nor parallelizable for finding Markovian Nash equilibrium, and a completely different approach is necessary.

\end{rem}

\subsection{Implementation}
A few details should be specified regarding the methodology described in Section \ref{sec_algorithm}.

First, the hypothesis spaces $\cN_0^{i'}$ and $\{\cN_k^i\}_{k=0}^{N_T-1}$ need to be specified. Note that, at each stage $m$, the optimal policy is defined through
$$\alpha^{i, m}(t, \bm x) = \argmin_{\alpha^i \in \mc{A}^i} G^i(t, \bm x, (\alpha^i, \bm \alpha^{-i,m-1}(t, \bm x)), \nabla_{\bx} V^{i, m}(t, \bx), V^{i, m}(t, \bx)),$$
one wishes to have direct access to the solved $V^{i,m}(t,\bx)$. Therefore we parametrize $V^{i}(t,\bx)$ (the superscript $m$ is dropped again for simplicity) directly with a neural network, denoted by $\text{Net}(t, \bx)$. Accordingly, $\text{Net}(0, \bx)$ becomes a hypothesis function in $\cN^{i'}_0$. For $\cN_k^i$, we know from the nonlinear Feynman-Kac formula \eqref{eq:FK relation} that $\Sigma(t,\bx)\nabla_{\bx} V^i(t,\bx)$ is the optimal map defining $\bZ_t^i$ in the variational problem \eqref{eq:variational_form}. Hence we choose $\Sigma(t_k,\bx)\nabla_{\bx} \text{Net}(t_k,\bx)$ to be a hypothesis function in $\cN_k^i$. In other words, the hypothesis functions in $\cN_0^{i'}$ and $\{\cN_k^i\}_{k=0}^{N_T-1}$ all share the same set of parameters. 
In this work, we use a fully-connected feedforward network with three hidden layers to instantiate $\text{Net}(t, \bx)$. The detailed architecture of $\text{Net}(t, \bx)$ are provided in Appendix~\ref{appendix_parameter}.

Second, at each stage $m$, we seek for an approximate solution to the PDE \eqref{eq:FP_PDE_explicit} for each player $i$, by iteratively updating the parameters of the neural networks. Note that the neural networks for $N$ players are decoupled and can be optimized in parallel. 
In our practice, the Adam method (\cite{Kingma2015adam}), a variant of SGD, is used to optimize the parameters. 
A relevant question is how many SGD iteration steps should be used per stage. After testing with different choices, we decide to use a moderate number as the iteration steps (see the discussion in Section~\ref{sec_example1} for details).
In fact, it is unwise to solve~\eqref{eq:FP_PDE_explicit} accurately at each stage with a lot of SGD updates, especially at early stages, as the opponents' strategies used for computing the best response are not even close to the Nash equilibrium. 
A similar idea of solving individual control problems not so accurately at each stage has also been used in \cite{SeBu:06}.
On the other hand, since the parameters are continued to be updated incrementally from the previous stage without re-initialization, it still suffices to expect the algorithm to converge with a moderate number of SGD updates per stage.

With the implementation details explained above, the pseudo-code of the proposed deep fictitious play algorithm is summarized in Algorithm \ref{def_algorithm1}.

\begin{algorithm}[!ht]
\caption{Deep Fictitious Play for Finding Markovian Nash Equilibrium \label{def_algorithm1}}
    \begin{algorithmic}[1]
	\REQUIRE $N$ = \# of players, $N_T$ = \# of subintervals on $[0,T]$, $M$ = \# of total stages in fictitious play, $N_{\text{sample}}$ = \# of sample paths generated for each player at each stage of fictitious play, $N_{\text{SGD\_per\_stage}}$ = \# of SGD steps for each player at each stage, $N_{\text{batch}}$ = batch size per SGD update, $\balpha^0\colon$ the initial policies that are smooth enough
	    \STATE Initialize $N$ deep neural networks to represent $V^{i,0}, i \in \mc{I}$
		\FOR{$m \gets 1$ to $M$}
		\FORALLP{$i \in \mc{I}$}
		\STATE  Generate $N_\text{sample}$ sample paths $\{\bX_{k}^{i,\pi}\}_{k=0}^{N_T}$ according to \eqref{eq:disc_X_path} and the realized optimal policies $\balpha^{-i, m-1}(t_k, \bm X_k^{i,\pi})$
		\FOR{$\ell \gets 1$ to $N_{\text{SGD}\_\text{per}\_\text{stage}}$}
		    \STATE Update the parameters of the $i^{th}$ neural network one step with $N_{\text{batch}}$ paths using the SGD algorithm (or its variant), based on the loss function \eqref{eq:disc_objective}
		  \ENDFOR
		  \STATE Obtain the approximate optimal policy  $\alpha^{i,m}$  according to \eqref{def_alphaast}
		  \ENDFOR
		  \STATE Collect the optimal policies at stage $m$: $\bm{\alpha}^m \gets (\alpha^{1,m}, \ldots, \alpha^{N,m})$
		  \ENDFOR
		\RETURN The optimal policy $\balpha^{M}$
	\end{algorithmic}
\end{algorithm}

\section{Numerical Examples}\label{sec_numerics}
In this section, we illustrate our algorithm on four examples, including games with identical or heterogeneous agents, and risk-neutral or risk-sensitive cost. In the main text we shall mainly focus on introducing the game setups and presenting numerical results. Technical details, such as concrete PDEs we aim to solve and the associated ground truth solutions $V^i$, will be kept to the minimal level. They, together with the hyperparameters and runtime of learning, are deferred to Appendices~\ref{appendix_PDE} and ~\ref{appendix_parameter}.

\subsection{An Inter-Bank Borrowing and Lending Game} 
\label{sec_example1}
Our first example models an inter-bank game concerning the systemic risk (\cite{CaFoSu:15}). Consider an inter-bank market with $N$ banks,
and let $X^i_t\in \RR$ be the log-monetary reserves of bank $i$ at time $t$. We model its dynamics as the following diffusion processes,
\begin{equation}
\label{eq:ex1_dynamics}
\ud X_t^i = [a(\overline X_t - X_t^i)   + \alpha_t^i ] \ud t  + \sigma \left(\rho \ud W_t^0 + \sqrt{1-\rho^2} W_t^i\right), \quad \overline X_t = \frac{1}{N}\sum_{i=1}^N X_t^i, \quad i \in \mc{I}.
\end{equation}
Here $a(\overline X_t - X_t^i)$ represents the rate at which bank $i$ borrows from or lends to other banks in the lending market, while $\alpha_t^i$  denotes its control rate of cash flows to a central bank. The standard Brownian motions $\{W_t^i\}_{i=0}^N$ are independent, in which $\{W_t^i, i\geq 1\}$ stands for the idiosyncratic noises and $W_t^0$ denotes the systemic shock, or so-called common noise in the general context. To describe the model in the form of \eqref{def_Xt}, we concatenate the log-monetary reserves $X_t^i$ of $N$ banks to form $\bX_t^{\balpha}=[X_t^1,\dots,X_t^N]\transpose$.
The associated drift term and diffusion term are defined as
\begin{equation}
\label{eq:example1_driftdiff}
b(t, \bm x, \balpha)=
[a(\bar x - x^1) + \alpha^1, 
\ldots, 
a(\bar x - x^N) + \alpha^N
]\transpose\in \RR^{N\times1}
, \quad \bar x = \frac{1}{N}\sum_{i=1}^N x^i, \end{equation}
\begin{equation}
\label{eq:example1_sigma}
\Sigma(t, \bm x)=
\begin{bmatrix}
     \sigma\rho & \sigma\sqrt{1-\rho^2} & 0 & \cdots & 0\\
     \sigma\rho & 0 & \sigma\sqrt{1-\rho^2} & \cdots & 0\\
     \vdots & \vdots & \vdots & \ddots & \vdots\\
     \sigma\rho & 0 & 0 & \cdots &\sigma\sqrt{1-\rho^2}
\end{bmatrix}\in \RR^{N\times (N+1)},
\end{equation}
and $\bW_t = (W_t^0, \ldots, W_t^N)$ is $(N+1)$-dimensional. The cost functional \eqref{def_J} that player $i$ wishes to minimize has the form
\begin{equation}
f^i(t, \bm{x},\balpha) = \half (\alpha^i)^2 - q \alpha^i(\bar x - x^i) + \frac{\eps}{2}(\bar x - x^i)^2,  \quad g^i(\bm{x}) = \frac{c}{2}(\bar x - x^i)^2.
\end{equation}
All the aforementioned parameters $a, \sigma, q, \eps, c$ are non-negative with $\abs{\rho} \leq 1$ and $q^2 \leq \eps$. We direct the interested readers to \cite{CaFoSu:15} for the detailed interpretation of this model.

\begin{figure}[!ht]
    \centering
    \includegraphics[width=\figwidth]{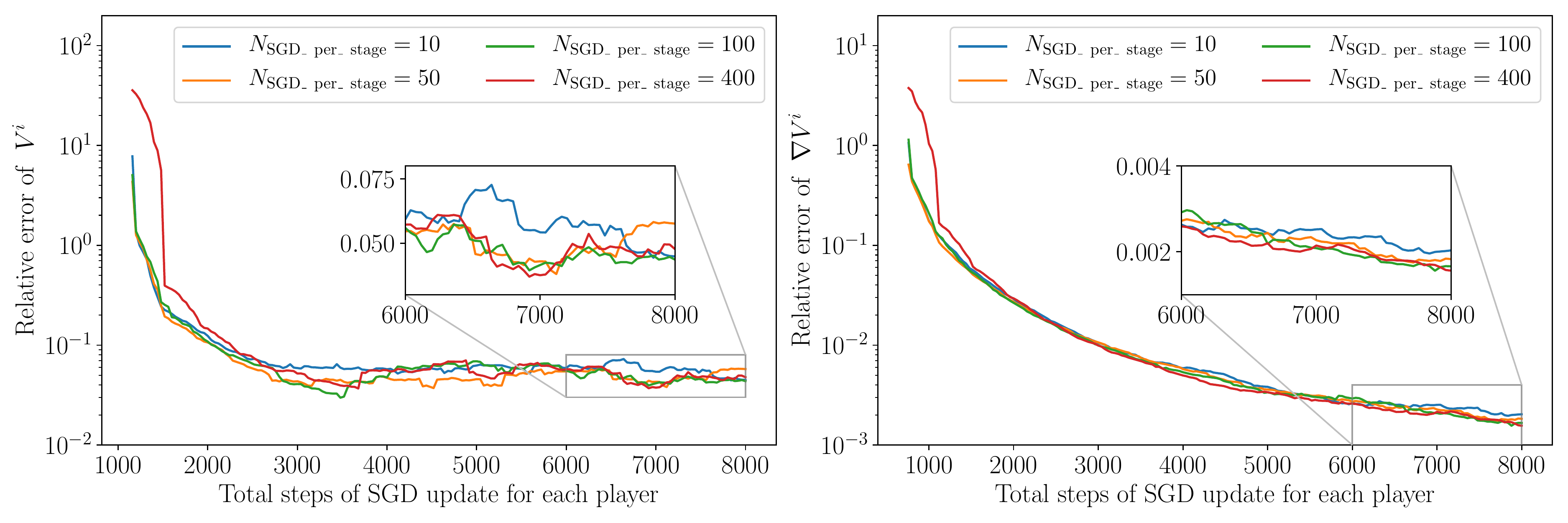}
    \caption{The relative squared errors of $V^i$ (left) and $\nabla V^i$ (right) along the training process of deep fictitious play for the inter-bank game in Section~\ref{sec_example1}. 
    The relative squared errors of $V^i(0, \bX_0^{i,\pi})$ and $\{\nabla V^{i}(t_k, \bX_k^{i,\pi})\}_{k=0}^{N_{T-1}}$ are evaluated.
    The error is computed every 400 SGD updates, averaged over all the players. A smoothed moving average with window size 3 is applied in the final plots.
    }
    \label{fig:LQMFG_err}
\end{figure}

The coupled HJB system corresponding to this game reads
\begin{align}\label{eq_HJBexample1}
\begin{dcases}
\partial_t V^i + \inf_{\alpha^i}\left\{\sum_{j=1}^N [a(\bar x - x^j) + \alpha^j] \partial_{x^j}V^i + \frac{(\alpha^i)^2}{2} - q\alpha^i(\bar x - x^i) + \frac{\eps}{2}(\bar x - x^i)^2 \right\} \\
\quad\quad\quad + \half \text{Tr}(\Sigma\transpose \text{Hess}_{\bx} V^i \Sigma) = 0,\\
V^i(T,\bm x) = \frac{c}{2}(\bar x - x^i)^2, \quad i \in \mc{I}.
\end{dcases}
\end{align} 
The minimizer in the infimum gives a candidate of the optimal control for player $i$:
$\alpha^{i}(t, \bm x) = q(\bar x - x^i) - \partial_{x^i}V^{i}(t, \bm x).$ 
Plugging it back into the $i^{th}$ equation yields a PDE of form \eqref{eq:FP_PDE_explicit}:
\begin{multline}\label{eq:example1}
    \partial_t V^{i} + \half \text{Tr}(\Sigma\transpose \text{Hess}_{\bx} V^i \Sigma) + a(\bar x - x^i)\partial_{x^i}V^{i} + \sum_{j\neq i} [a(\bar x - x^j) + \alpha^{j}(t,\bm x)] \partial_{x^j}V^{i} \\
 + \frac{\eps}{2}(\bar x - x^i)^2 - \half(q(\bar x - x^i) - \partial_{x^i} V^{i})^2= 0,
\end{multline}
where $\alpha^j$ with $j\neq i$ are considered exogenous for player $i$'s problem, and are given by the best responses of the other players from the previous stage. To be precise, $\mu^i$ and $h^i$ in \eqref{eq:FP_PDE_explicit} are defined as:
\begin{align}
\mu^i(t, \bx; \balpha^{-i})&=
[
a(\bar{x} - x^1) + \alpha^1, 
\ldots,
a(\bar{x} - x^i),
\ldots,
a(\bar{x} - x^N) + \alpha^N
]\transpose,\\
h^i(t, \bx,  y, \bz;  \balpha^{-i})&=\frac{\eps}{2}(\bar x - x^i)^2 - \half(q(\bar x - x^i) - \frac{z^i}{\sigma\sqrt{1-\rho^2}})^2,
\end{align}
in which $\bz=(z^0, z^1,\dots, z^N)\in\RR^{N+1}$. 
\begin{figure}[!ht]
    \centering
    \includegraphics[width=\figwidth]{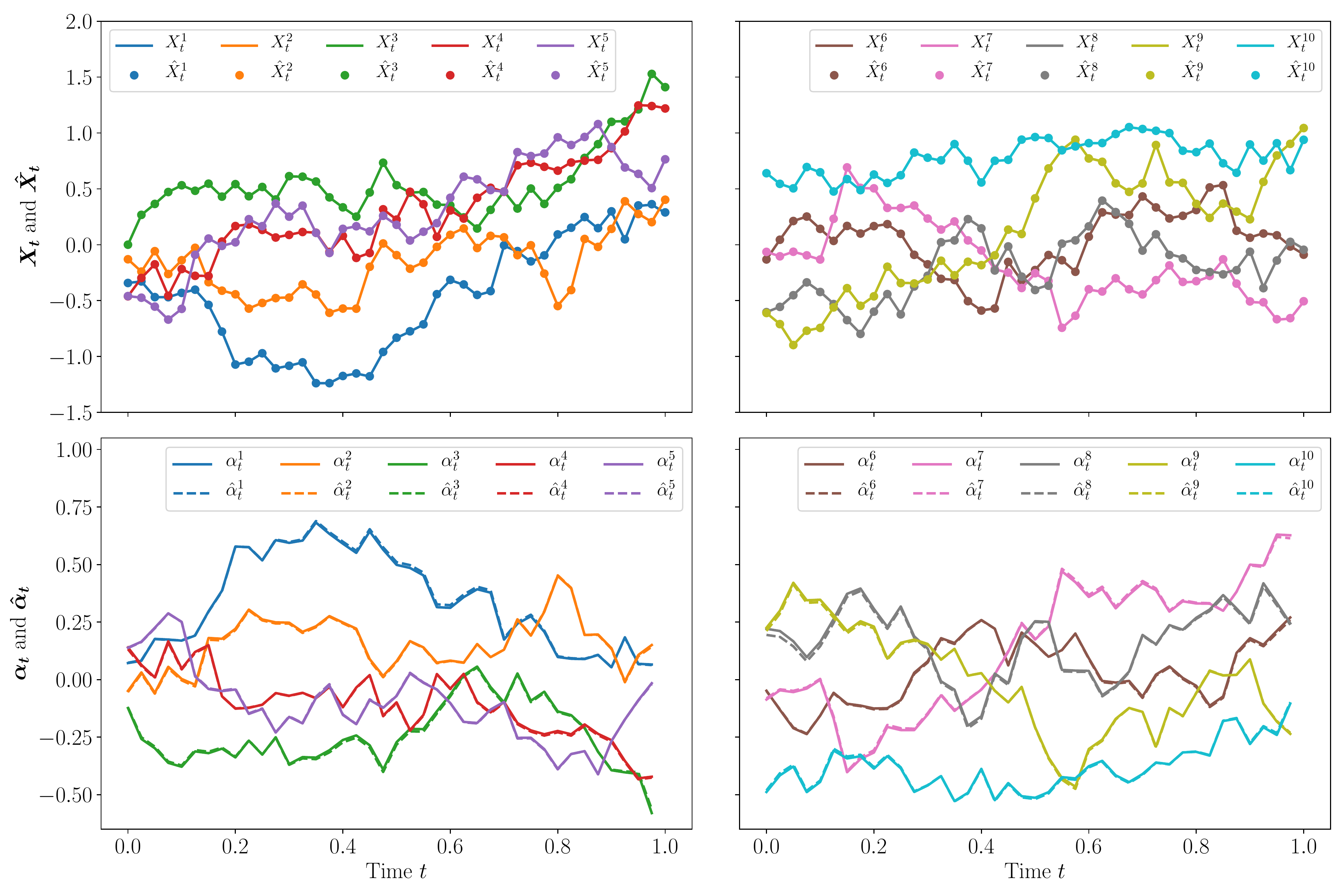}
    \caption{A sample path for each player of the inter-bank game in Section~\ref{sec_example1} with $N=10$. Top: the optimal state process $X_t^i$ (solid lines) and its approximation $\hat{X}_t^i$ (circles) provided by the optimized neural networks, under the same realized path of Brownian motion. Bottom: comparisons of the strategies $\alpha_t^i$ and $\hat{\alpha}_t^i$ (dashed lines).}
    \label{fig:LQMFG_path}
\end{figure}

Figures~\ref{fig:LQMFG_err}--\ref{fig:LQMFG_path} show the performance of our algorithm on a 10-player game, using the parameter:
\begin{equation}\label{def_parameters}
    a = 0.1,\quad  q = 0.1, \quad c = 0.5,\quad  \eps = 0.5,\quad  \rho = 0.2, \quad \sigma = 1, \quad T = 1.
\end{equation}
We define the relative squared error by
{\footnotesize
\begin{align*}
    \text{RSE} = \frac{\sum_{\substack{i \in \mathcal{I}\\1 \leq  j \leq J}} \left(V^i(0, \bm x_{t_0}^{(j)}) - \widehat V^i(0,\bm x_{t_0}^{(j)})\right)^2}{\sum_{\substack{i \in \mathcal{I}\\ 1 \leq j \leq J}} \left(V^i(0, \bm x_{t_0}^{(j)}) - \bar V^i\right)^2},
    \;\text{or }\; 
    \text{RSE} = \frac{\sum_{\substack{i \in \mathcal{I}\\ 0 \leq k \leq N_T-1 \\ 1 \leq j \leq J}} \left(\nabla_{\bx} V^i(t_k, \bm x_{t_k}^{(j)}) - \nabla_{\bx} \widehat V^i(t_k,\bm x_{t_k}^{(j)})\right)^2}{\sum_{\substack{i \in \mathcal{I}\\ 0 \leq k \leq N_T-1\\1 \leq j \leq J}} \left(\nabla_{\bx} V^i(t_k, \bm x_{t_k}^{(j)})- \overline {\nabla_{\bx} V}^i\right)^2},
\end{align*}
}where 
$V^i$ is given by the explicit formula provided in Appendix~\ref{appendix_example1},
$\hat{V}^i$ is the prediction from the neural networks,
and $\bar V^i ~(\emph{resp.~} \overline {\nabla_{\bx} V}^i)$ is the average of $V^i ~(\emph{resp.~} \nabla_{\bx}V^i)$ evaluated at all the indices $j, k$.
To compute the relative error, we generate $J=256$ ground truth sample paths $\{\bm x_{t_k}^{(j)}\}_{k=0}^{N_T-1}$ using Euler scheme based on  \eqref{def_Xt}\eqref{eq:example1_driftdiff}\eqref{eq:example1_sigma} and the true optimal strategy (provided in Appendix~\ref{appendix_example1}). Note that the superscript ${(j)}$ here does not mean the player index, but the $j^{th}$ path for all players. 
The relative errors reported in Sections~\ref{sec_example2} and \ref{sec_example3} are defined in the same way.
Figure~\ref{fig:LQMFG_err} in particular compares the relative squared error
as $N_{\text{SGD\_per\_stage}}$ varies from 10 to 400. The convergence of the learning curves with small $N_{\text{SGD\_per\_stage}}$ asserts that each individual problem does not need to be solved so accurately. Furthermore, the similar performances under different $N_{\text{SGD\_per\_stage}}$ with the same total budget of SGD updates suggests that the algorithm is insensitive to the choice of this hyperparameter. In all the following numerical experiments, we fix $N_{\text{SGD\_per\_stage}}=100$.
The final relative squared errors of $V$ and $\nabla V$ averaged from three independent runs of deep fictitious play are 4.6\% and 0.2\%,  respectively. 
Figure~\ref{fig:LQMFG_path} presents one sample path for each player of the optimal state process $X_t^i$ and the optimal control $\alpha_t^i$ \emph{vs.} their approximations $\hat{X}_t^i, \hat{\alpha}_t^i$ provided by the optimized neural networks.

One concern is that how sensitive the numerical result is to the parameters chosen in \eqref{def_parameters}. For instance, in solving mean-field games, the numerical algorithms may produce  bifurcations (cf. \cite{AnCh:18,ChCrDe:19}). To show the robustness (within a certain range of parameters) of our algorithm, we conduct two experiments by increasing the game length $T$, and the coupling between players described by $a$. The results are provided in Appendix~\ref{appendix_sensitivity}.


\subsection{A Risk-Sensitive Version of the Inter-Bank Game in Section~\ref{sec_example1}}
\label{sec_example2}

Next we consider the same inter-bank game as in Section~\ref{sec_example1}, but under risk-sensitive utility \eqref{def_J'}. The cases $\theta_i <0$, $\theta_i >0$ correspond respectively to what an economist would term \emph{risk-preferring} and \emph{risk-averse} attitudes on the expected utility. Finding the Nash equilibrium of this game is reduced to solving $N$-coupled matrix Riccati equations. Its derivation, as well as the concrete form of PDE \eqref{eq:FP_PDE_explicit} is presented in Appendix~\ref{appendix_example2}.

\begin{figure}[!ht]
    \centering
    \includegraphics[width=\figwidth]{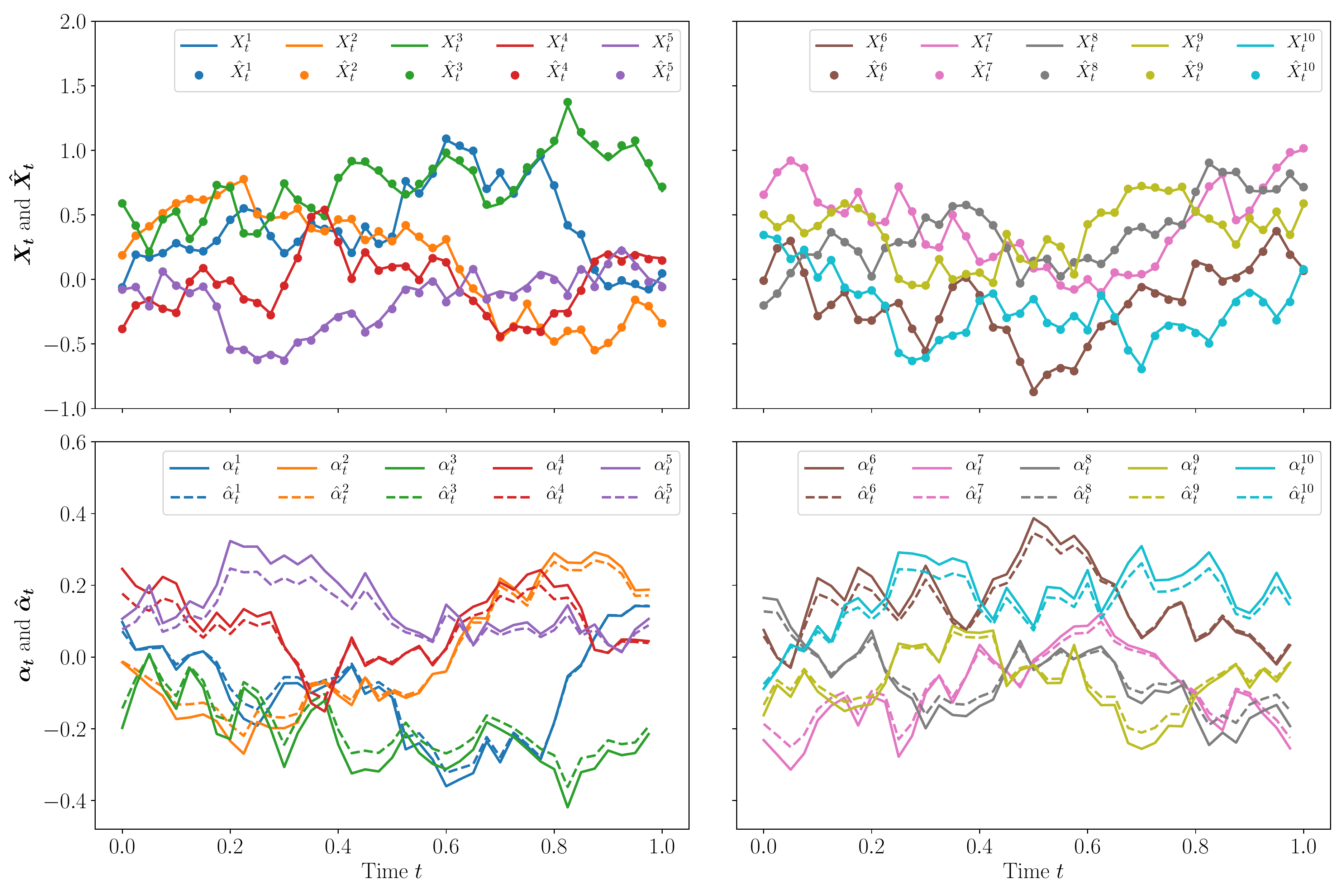}
    \caption{A sample path for each player of the risk-sensitive inter-bank game in Section~\ref{sec_example2} with $N=10$. Top: the optimal state process $X_t^i$ (solid lines) and its approximation $\hat{X}_t^i$ (circles) provided by the optimized neural networks, under the same realized path of Brownian motion. Bottom: comparisons of the strategies $\alpha_t^i$ and $\hat{\alpha}_t^i$ (dashed lines).}
    \label{fig:RSLQMFG_path}
\end{figure}

For numerical illustration, we study a game of 10 heterogeneous players, with risk-sensitivity $\theta_i = 0.6 + 0.02i$. Other parameters follow \eqref{def_parameters} except $c=0.3, \eps=0.3$. As before, Figure~\ref{fig:RSLQMFG_path} compares the sample paths of the true optimal states/control \emph{vs.} the approximated ones produced by the proposed deep fictitious play algorithm. 
The final relative squared errors of $V$ and $\nabla V$ averaged from three independent runs are 1.4\% and 0.1\%, respectively.

\subsection{A General Linear-Quadratic Risk-Sensitive Dynamic Game}
\label{sec_example3}
In this section, we consider a general linear-exponential-quadratic game of the form \eqref{def_Xt}:
\begin{equation}
b(t, \bm x, \bm \alpha) = A\bm x + \sum_{i=1}^N B_i \alpha^i(t, \bm x), \quad \Sigma(t, \bm x) = \Sigma,
\end{equation}
where $A \in \RR^{n\times n}$, $B_i \in \RR^{n \times n_i}$, $\Sigma \in \RR^{n \times n}$ are constant matrices, $\Sigma$ is of full rank, $\bx \in \RR^n, \alpha^i\colon [0,T] \times \RR^n \to \RR^{n_i}$. 
In contrast to the first two examples in which each player $i$ can only control her own state $X^i$, here the control $\alpha^i$ of each player $i$ contributes to the dynamics of the common state $\bX$ through a general linear relationship. The cost functional for each player is risk-sensitive of the form \eqref{def_J'} with:
\begin{equation}
f^i(t, \bm x, \bm \alpha) = \half\bm x\transpose Q_i \bm x + \half(\alpha^i)\transpose R_i \alpha^i, \quad g^i(\bm x) = \half\bm x\transpose M_i \bm x,     
\end{equation}
where $Q_i \in \RR^{n \times n}$, $R_i \in \RR^{n_i \times n_i}$ and $M_i \in \RR^{n \times n}$ are symmetric, positive definite constant matrices. 

\begin{figure}[ht]
    \centering
    \includegraphics[width=\figwidth]{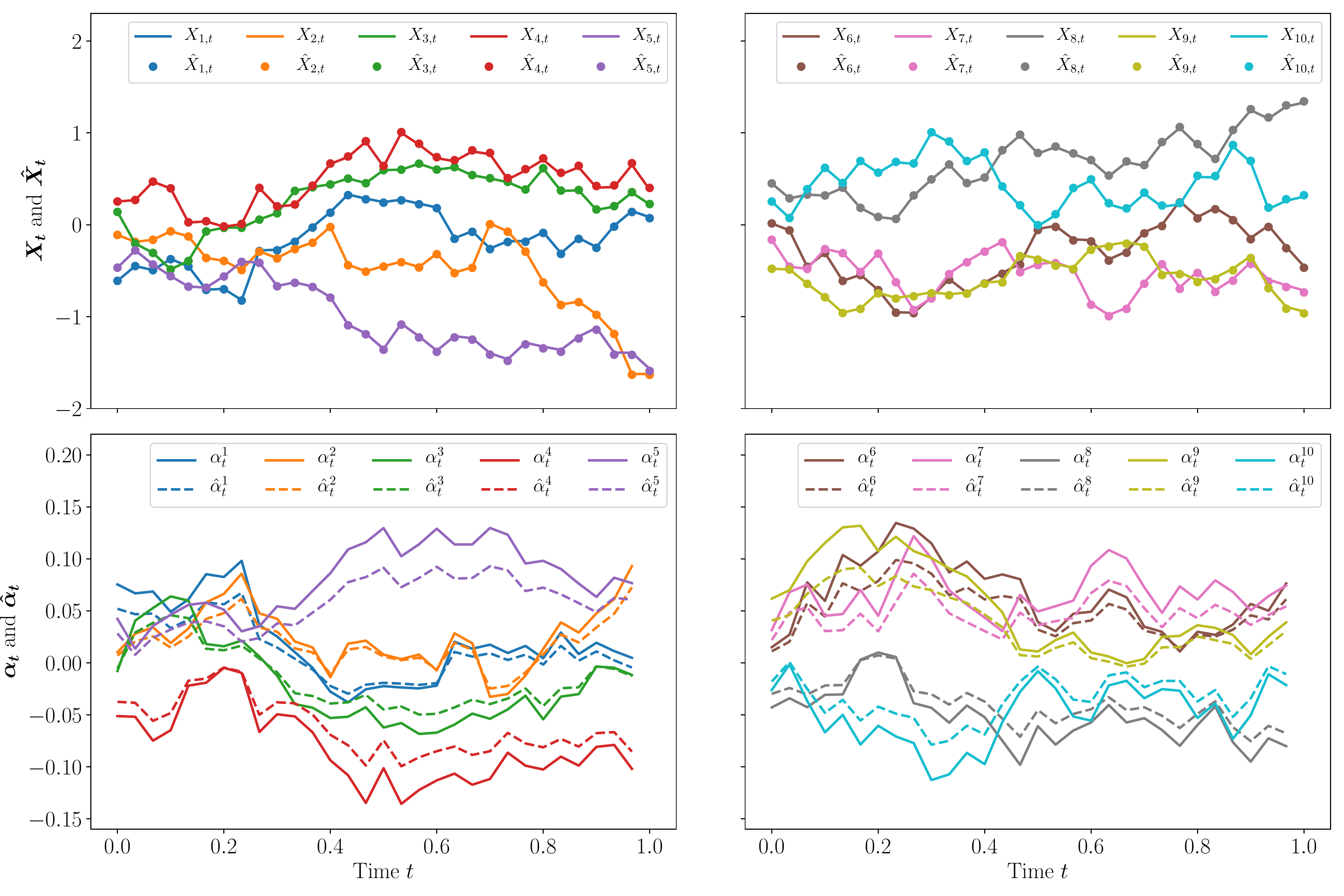}
    \caption{A sample path for each player of the general linear-exponential-quadratic game in Section~\ref{sec_example3} with $N=10$. Top: the optimal state process $X_{i,t}$ (solid lines) and its approximation $\hat{X}_{i,t}$ (circles) provided by the optimized neural networks, under the same realized path of Brownian motion. Bottom: comparisons of the strategies $\alpha_t^i$ and $\hat{\alpha}_t^i$ (dashed lines).}
    \label{fig:RSGMat_path}
\end{figure}

As before, one can solve $N$-coupled matrix Riccati equations to obtain the Nash equilibrium. The  technical details are deferred to Appendix~\ref{appendix_example3}, and we shall focus on the numerics. We present a 10-player heterogeneous game in which the dimensions of the state $\bX_t$ and the controls $\alpha^i$ are all 10, i.e., $n=n_i=10$. The risk-sensitive parameters are the same as in Section \ref{sec_example2}, i.e.,  $\theta_i=0.6+0.02i$. Denote by $I_{10}\in\RR^{10\times 10}$ the identity matrix, and let $B_i=R_i=0.4I_{10}, \Sigma=I_{10}$. The matrices $A, Q_i, M_i$ are defined by
\begin{align*}
    A = 0.1I_{10} + 0.05(\tilde{A} + \tilde{A}\transpose), \quad
    Q_i = 0.1I_{10} + 0.05(\tilde{Q_i} + \tilde{Q_i}\transpose), \quad
    M_i = 0.2I_{10} + 0.1(\tilde{M_i} + \tilde{M_i}\transpose).
\end{align*}
Here $\tilde{A}, \tilde{Q_i}, \tilde{M_i}$ are all random matrices whose diagonal entries are 0 and off-diagonal entries are independently sampled from the uniform distribution on $[-1, 1]$.
Hence the game is entirely heterogeneous due to different $\theta_i, Q_i, M_i$.
As before, we illustrate the sample paths of the true optimal state process/policy \emph{vs.} the approximated ones in Figure~\ref{fig:RSGMat_path}. The final relative squared errors of $V$ and $\nabla V$ averaged from three independent trials are 6.5\% and 0.4\%, respectively. 

\subsection{A Variant of the Inter-Bank Game in Section~\ref{sec_example1} with $ \text{\boldmath$N=50$}$}
\label{sec_example4}
\begin{figure}[ht]
    \centering
    \includegraphics[width=\figwidth]{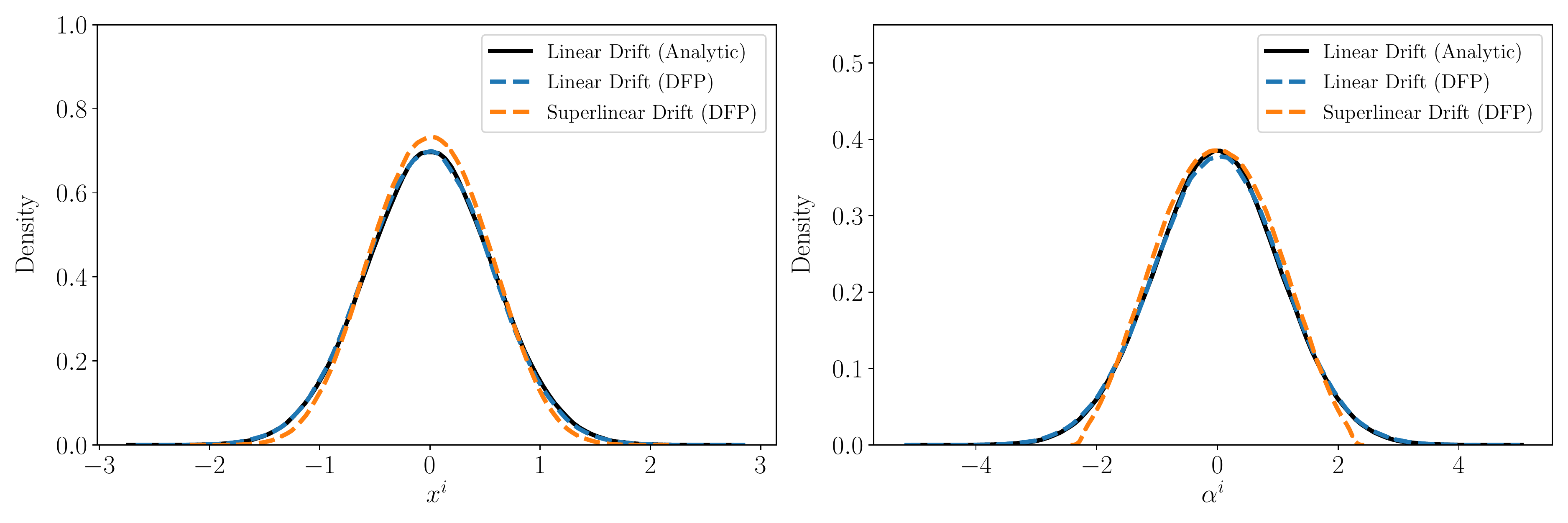}
    \caption{The densities of $X^i_T$ (left) and $\alpha^i_T$ (right) of the inter-bank game with $N=50$. We plot the true distribution of the linear dynamics~\eqref{eq:ex1_dynamics} (black solid lines), and the distributions obtained from the deep fictitious play with the linear dynamics~\eqref{eq:ex1_dynamics} (blue dashed lines) and superlinear dynamics ~\eqref{eq:ex4_dynamics} (orange dashed lines).}
    \label{fig:NonLQMFG_hist}
\end{figure}

In this section we first recompute the example in Section~\ref{sec_example1} with all the same parameters as in \eqref{def_parameters} but $N=50$ to test the performance of the proposed deep fictitious play with even more players. 
Since this game is totally symmetric for all the players, the fictitious play in fact reduces it into $N$ identical decision problems. Due to this, at each stage of the fictitious play, if all the players use the same strategy, then the updated strategies for them are still essentially the same. Therefore we can leverage such symmetry of the game to reduce the computational cost when $N$ is large. Specifically, we impose that all the players always share the same strategy at each stage
by using a single network to derive each player's strategy. Consequently, at each stage, only one player's decision problem needs to be solved by the deep BSDE method to update the neural network-based strategy. 
Algorithm~\ref{def_algorithm1} is accordingly simplified, whose details are provided in the Appendix~\ref{appendix_parameter}. After solving  \eqref{eq_HJBexample1} with the deep fictitious play, we simulate 10,000 paths following the optimized neural network, and we plot the histograms of $\hat \bX_T$ and $\hat \balpha_T$ in Figure~\ref{fig:NonLQMFG_hist} by considering all 50 components together, i.e., the histograms are generated by 500,000 data points. Note that due to all the symmetry of this problem, the distribution of each component is identical. Therefore considering them together empirically does not change the underlying distribution but reduces the variance. Figure~\ref{fig:NonLQMFG_hist} shows great consistency between the histograms obtained from the optimized strategies/paths (blue dashed lines) and true optimal strategies/paths (black solid lines).

Convinced by the reliability of the algorithm with a large number of players, we further consider a variant of the example in Section~\ref{sec_example1} with $N=50$ again.
In contrast to \eqref{eq:ex1_dynamics},
the drift term becomes nonlinear
\begin{equation}
\label{eq:ex4_dynamics}
\ud X_t^i = [a(\overline X_t - X_t^i)^3  + \alpha_t^i ] \ud t  + \sigma \left(\rho \ud W_t^0 + \sqrt{1-\rho^2} W_t^i\right).
\end{equation}
In this case, as far as the authors are aware, there is no analytic solution or simple characterization of the Nash equilibrium suitable for numerical computation. For the corresponding mean-field game, due to the presence of common noise, the Nash equilibrium is characterized by a coupled system of stochastic partial differential equations, which is also very difficult, if not impossible, to solve numerically. In computation of the deep fictitious play, we set $a=10$ in order to compensate the smaller drift near 0 and all other parameters the same as in \eqref{def_parameters}. It turns out that the algorithm still finds reasonable equilibrium for this problem, as shown in Figure~\ref{fig:NonLQMFG_hist} as well. For $X_T^i$, since the drift term $a(\overline X_t - X_t^i)^3$ is superlinear, the final distribution of $X_T^i$ is expected to be more concentrated than the one under linear drift, which is normal distributed with kurtosis 3. It is confirmed in the left panel of Figure~\ref{fig:NonLQMFG_hist}, in which the kurtosis of the orange dashed line is $2.72 < 3$. 
For $\alpha^i_T$, if $X^i_t$ is far away from the average, the superlinear term $a(\overline X_t - X_t^i)^3$ will push it back quickly, saving some effort of bank $i$ and reducing $\alpha^i_t$. Therefore 
the tail of $\alpha^i_T$ is much lighter than the Gaussian distribution (with kurtosis $2.36 < 3$), as shown in the right panel of Figure~\ref{fig:NonLQMFG_hist}. 

\section{Conclusion}\label{sec_conclusion}
In this paper, we propose a deep fictitious play algorithm to compute the Markovian Nash equilibrium of large $N$-player stochastic differential games. The game is firstly decoupled into $N$ individual decision problems by the idea of fictitious play, and then each is solved iteratively.
Due to the feedback nature of Markovian Nash equilibrium, it is inefficient to directly parameterize the optimal policy as in \cite{Hu2:19}.
Instead, we rely on the HJB approach and solve it through the deep BSDE method.
Three examples with closed-form solutions are carefully examined, and the algorithm performs unanimously well. The fourth example without tractability is also presented evidencing the algorithm's applicability to general games without the linear-(exponential)-quadratic structure.

Future work includes the convergence analysis on the proposed algorithm and more practical applications in other disciplines such as operation research and economics. For stochastic differential games involving controlling volatility, it is also promising to use the deep fictitious play to find the Nash equilibrium, by combining a second-order BSDE formulation.

\acks{Part of this work was done during the visit of JH and RH to the Beijing Institute of Big Data
Research, China. They really appreciate the hospitality of the institute, and thank Professor Weinan
E for hosting and useful discussions.}

\bibliography{Reference}

\begin{thebibliography}{56}
\providecommand{\natexlab}[1]{#1}
\providecommand{\url}[1]{\texttt{#1}}
\expandafter\ifx\csname urlstyle\endcsname\relax
  \providecommand{\doi}[1]{doi: #1}\else
  \providecommand{\doi}{doi: \begingroup \urlstyle{rm}\Url}\fi

\bibitem[Achdou et~al.(2017)Achdou, Han, Lasry, Lions, and Moll]{HACT2017}
Y.~Achdou, J.~Han, J.-M. Lasry, P.-L. Lions, and B.~Moll.
\newblock Income and wealth distribution in macroeconomics: A continuous-time
  approach.
\newblock Working Paper 23732, National Bureau of Economic Research, August
  2017.

\bibitem[Angiuli et~al.(2018)Angiuli, Graves, Li, Chassagneux, Delarue, and
  Carmona]{AnCh:18}
A.~Angiuli, C.~V. Graves, H.~Li, J.-F. Chassagneux, F.~Delarue, and R.~Carmona.
\newblock Numerical probabilistic approach to mfg.
\newblock \emph{arXiv preprint arXiv:1805.02406}, 2018.

\bibitem[Bengio(2009)]{Be:09}
Y.~Bengio.
\newblock Learning deep architectures for {AI}.
\newblock \emph{Foundations and trends{\textregistered} in Machine Learning},
  2\penalty0 (1):\penalty0 1--127, 2009.

\bibitem[Bensoussan et~al.(2014)Bensoussan, Siu, Yam, and Yang]{BeSiYaYa:14}
A.~Bensoussan, C.~C. Siu, S.~C.~P. Yam, and H.~Yang.
\newblock A class of non-zero-sum stochastic differential investment and
  reinsurance games.
\newblock \emph{Automatica}, 50\penalty0 (8):\penalty0 2025--2037, 2014.

\bibitem[Berger(2005)]{Be:05}
U.~Berger.
\newblock Fictitious play in 2$\times$ n games.
\newblock \emph{Journal of Economic Theory}, 120\penalty0 (2):\penalty0
  139--154, 2005.

\bibitem[Briani and Cardaliaguet(2018)]{BrCa:18}
A.~Briani and P.~Cardaliaguet.
\newblock Stable solutions in potential mean field game systems.
\newblock \emph{Nonlinear Differential Equations and Applications}, 25\penalty0
  (1):\penalty0 1, 2018.

\bibitem[Brown(1949)]{Br:49}
G.~W. Brown.
\newblock Some notes on computation of games solutions.
\newblock Technical report, Rand Corp Santa Monica CA, 1949.

\bibitem[Brown(1951)]{Br:51}
G.~W. Brown.
\newblock Iterative solution of games by fictitious play.
\newblock \emph{Activity Analysis of Production and Allocation}, 13\penalty0
  (1):\penalty0 374--376, 1951.

\bibitem[{Busoniu} et~al.(2008){Busoniu}, {Babuska}, and {De
  Schutter}]{Busoniu08}
L.~{Busoniu}, R.~{Babuska}, and B.~{De Schutter}.
\newblock A comprehensive survey of multiagent reinforcement learning.
\newblock \emph{IEEE Transactions on Systems, Man, and Cybernetics, Part C
  (Applications and Reviews)}, 38\penalty0 (2):\penalty0 156--172, 2008.

\bibitem[Cardaliaguet and Hadikhanloo(2017)]{CaHa:17}
P.~Cardaliaguet and S.~Hadikhanloo.
\newblock Learning in mean field games: the fictitious play.
\newblock \emph{ESAIM: Control, Optimisation and Calculus of Variations},
  23\penalty0 (2):\penalty0 569--591, 2017.

\bibitem[Cardaliaguet and Lehalle(2018)]{CaLe:18}
P.~Cardaliaguet and C.-A. Lehalle.
\newblock Mean field game of controls and an application to trade crowding.
\newblock \emph{Mathematics and Financial Economics}, 12\penalty0 (3):\penalty0
  335--363, 2018.

\bibitem[Carleo and Troyer(2017)]{CaTr:17}
G.~Carleo and M.~Troyer.
\newblock Solving the quantum many-body problem with artificial neural
  networks.
\newblock \emph{Science}, 355\penalty0 (6325):\penalty0 602--606, 2017.

\bibitem[Carmona and Delarue(2013)]{CaDe:13}
R.~Carmona and F.~Delarue.
\newblock Probabilistic analysis of mean-field games.
\newblock \emph{SIAM Journal on Control and Optimization}, 51\penalty0
  (4):\penalty0 2705--2734, 2013.

\bibitem[Carmona et~al.(2015)Carmona, Fouque, and Sun]{CaFoSu:15}
R.~Carmona, J.-P. Fouque, and L.-H. Sun.
\newblock Mean field games and systemic risk.
\newblock \emph{Communications in Mathematical Sciences}, 13\penalty0
  (4):\penalty0 911--933, 2015.

\bibitem[Chassagneux et~al.(2019)Chassagneux, Crisan, and Delarue]{ChCrDe:19}
J.-F. Chassagneux, D.~Crisan, and F.~Delarue.
\newblock Numerical method for fbsdes of mckean--vlasov type.
\newblock \emph{The Annals of Applied Probability}, 29\penalty0 (3):\penalty0
  1640--1684, 2019.

\bibitem[Dockner et~al.(2000)Dockner, Jorgensen, Van~Long, and
  Sorger]{DoJoVaSo:00}
E.~J. Dockner, S.~Jorgensen, N.~Van~Long, and G.~Sorger.
\newblock \emph{Differential Games in Economics and Management Science}.
\newblock Cambridge University Press, 2000.

\bibitem[E et~al.(2017)E, Han, and Jentzen]{EHaJe:17}
W.~E, J.~Han, and A.~Jentzen.
\newblock Deep learning-based numerical methods for high-dimensional parabolic
  partial differential equations and backward stochastic differential
  equations.
\newblock \emph{Communications in Mathematics and Statistics}, 5\penalty0
  (4):\penalty0 349--380, 2017.

\bibitem[El~Karoui et~al.(1997)El~Karoui, Peng, and Quenez]{ElPeQu:97}
N.~El~Karoui, S.~Peng, and M.~C. Quenez.
\newblock Backward stochastic differential equations in finance.
\newblock \emph{Mathematical Finance}, 7\penalty0 (1):\penalty0 1--71, 1997.

\bibitem[Foerster et~al.(2016)Foerster, Assael, de~Freitas, and
  Whiteson]{foerster2016learning}
J.~Foerster, I.~A. Assael, N.~de~Freitas, and S.~Whiteson.
\newblock Learning to communicate with deep multi-agent reinforcement learning.
\newblock In \emph{Advances in Neural Information Processing Systems}, pages
  2137--2145, 2016.

\bibitem[Foster and Young(1998)]{FoYo:98}
D.~P. Foster and H.~P. Young.
\newblock On the nonconvergence of fictitious play in coordination games.
\newblock \emph{Games and Economic Behavior}, 25\penalty0 (1):\penalty0 79--96,
  1998.

\bibitem[Gao and Duan(2017)]{GaDu:17}
X.~Gao and L.-M. Duan.
\newblock Efficient representation of quantum many-body states with deep neural
  networks.
\newblock \emph{Nature Communications}, 8\penalty0 (1):\penalty0 662, 2017.

\bibitem[Gomes and Voskanyan(2013)]{GoVo:13}
D.~A. Gomes and V.~K. Voskanyan.
\newblock Extended deterministic mean-field games.
\newblock \emph{arXiv preprint arXiv:1305.2600}, 2013.

\bibitem[Gomes and Voskanyan(2016)]{GoVo:16}
D.~A. Gomes and V.~K. Voskanyan.
\newblock Extended deterministic mean-field games.
\newblock \emph{SIAM Journal on Control and Optimization}, 54\penalty0
  (2):\penalty0 1030--1055, 2016.

\bibitem[Gomes et~al.(2014)Gomes, Patrizi, and Voskanyan]{GoPaVo:14}
D.~A. Gomes, S.~Patrizi, and V.~Voskanyan.
\newblock On the existence of classical solutions for stationary extended mean
  field games.
\newblock \emph{Nonlinear Analysis: Theory, Methods \& Applications},
  99:\penalty0 49--79, 2014.

\bibitem[Han and Long(2018)]{HaLo:18}
J.~Han and J.~Long.
\newblock Convergence of the deep {BSDE} method for coupled {FBSDE}s.
\newblock \emph{arXiv:1811.01165}, 2018.

\bibitem[Han et~al.(2018)Han, Jentzen, and E]{HaJeE:18}
J.~Han, A.~Jentzen, and W.~E.
\newblock Solving high-dimensional partial differential equations using deep
  learning.
\newblock \emph{Proceedings of the National Academy of Sciences}, 115\penalty0
  (34):\penalty0 8505--8510, 2018.

\bibitem[Han et~al.(2019)Han, Ma, Ma, and E]{han2019uniformly}
J.~Han, C.~Ma, Z.~Ma, and W.~E.
\newblock Uniformly accurate machine learning-based hydrodynamic models for
  kinetic equations.
\newblock \emph{Proceedings of the National Academy of Sciences}, 116\penalty0
  (44):\penalty0 21983--21991, 2019.

\bibitem[Hofbauer and Sandholm(2002)]{HoSa:02}
J.~Hofbauer and W.~H. Sandholm.
\newblock On the global convergence of stochastic fictitious play.
\newblock \emph{Econometrica}, 70\penalty0 (6):\penalty0 2265--2294, 2002.

\bibitem[Hu(2019)]{Hu2:19}
R.~Hu.
\newblock Deep fictitious play for stochastic differential games.
\newblock \emph{arXiv preprint arXiv:1903.09376}, 2019.

\bibitem[Huang et~al.(2006)Huang, Malham{\'e}, and Caines]{HuMaCa:06}
M.~Huang, R.~P. Malham{\'e}, and P.~E. Caines.
\newblock Large population stochastic dynamic games: closed-loop
  {McKean-Vlasov} systems and the {Nash} certainty equivalence principle.
\newblock \emph{Communications in Information and Systems}, 6\penalty0
  (3):\penalty0 221--252, 2006.

\bibitem[Huang et~al.(2007)Huang, Caines, and Malham{\'e}]{HuCaMa:07}
M.~Huang, P.~E. Caines, and R.~P. Malham{\'e}.
\newblock Large-population cost-coupled {LQG} problems with nonuniform agents:
  individual-mass behavior and decentralized $\epsilon$-{Nash} equilibria.
\newblock \emph{IEEE Transactions on Automatic Control}, 52\penalty0
  (9):\penalty0 1560--1571, 2007.

\bibitem[Ioffe and Szegedy(2015)]{batch}
S.~Ioffe and C.~Szegedy.
\newblock Batch normalization: Accelerating deep network training by reducing
  internal covariate shift.
\newblock In \emph{International Conference on Machine Learning}, pages
  448--456, 2015.

\bibitem[Jordan(1993)]{Jo:93}
J.~S. Jordan.
\newblock Three problems in learning mixed-strategy {Nash} equilibria.
\newblock \emph{Games and Economic Behavior}, 5\penalty0 (3):\penalty0
  368--386, 1993.

\bibitem[Kingma and Ba(2015)]{Kingma2015adam}
D.~Kingma and J.~Ba.
\newblock Adam: a method for stochastic optimization.
\newblock In \emph{Proceedings of the International Conference on Learning
  Representations}, 2015.

\bibitem[Krishna and Sj{\"o}str{\"o}m(1998)]{KrSj:98}
V.~Krishna and T.~Sj{\"o}str{\"o}m.
\newblock On the convergence of fictitious play.
\newblock \emph{Mathematics of Operations Research}, 23\penalty0 (2):\penalty0
  479--511, 1998.

\bibitem[Lasry and Lions(2006{\natexlab{a}})]{LaLi1:2006}
J.-M. Lasry and P.-L. Lions.
\newblock Jeux à champ moyen. {I.} {L}e cas stationnaire.
\newblock \emph{C. R. Math. Acad. Sci. Paris}, 9:\penalty0 619--625,
  2006{\natexlab{a}}.

\bibitem[Lasry and Lions(2006{\natexlab{b}})]{LaLi2:2006}
J.-M. Lasry and P.-L. Lions.
\newblock Jeux à champ moyen. {II.} {H}orizon fini et contrôle optimal.
\newblock \emph{C. R. Math. Acad. Sci. Paris}, 10:\penalty0 679--684,
  2006{\natexlab{b}}.

\bibitem[Lasry and Lions(2007)]{LaLi:2007}
J.-M. Lasry and P.-L. Lions.
\newblock Mean field games.
\newblock \emph{Japanese Journal of Mathematics}, 2:\penalty0 229--260, 2007.

\bibitem[LeCun et~al.(2015)LeCun, Bengio, and Hinton]{LeBeHi:15}
Y.~LeCun, Y.~Bengio, and G.~Hinton.
\newblock Deep learning.
\newblock \emph{Nature}, 521\penalty0 (7553):\penalty0 436, 2015.

\bibitem[Milgrom and Roberts(1991)]{MiRo:91}
P.~Milgrom and J.~Roberts.
\newblock Adaptive and sophisticated learning in normal form games.
\newblock \emph{Games and Economic Behavior}, 3\penalty0 (1):\penalty0 82--100,
  1991.

\bibitem[Miyasawa(1961)]{Mi:61}
K.~Miyasawa.
\newblock On the convergence of the learning process in a $2\times 2$
  non-zero-sum two-person game.
\newblock Technical report, Princeton University NJ, 1961.

\bibitem[Monderer and Sela(1996)]{MoSe:96}
D.~Monderer and A.~Sela.
\newblock A 2$\times$ 2 game without the fictitious play property.
\newblock \emph{Games and Economic Behavior}, 14\penalty0 (1):\penalty0
  144--148, 1996.

\bibitem[Monderer and Shapley(1996{\natexlab{a}})]{MoSh2:96}
D.~Monderer and L.~S. Shapley.
\newblock Potential games.
\newblock \emph{Games and Economic Behavior}, 14\penalty0 (1):\penalty0
  124--143, 1996{\natexlab{a}}.

\bibitem[Monderer and Shapley(1996{\natexlab{b}})]{MoSh:96}
D.~Monderer and L.~S. Shapley.
\newblock Fictitious play property for games with identical interests.
\newblock \emph{Journal of Economic Theory}, 68\penalty0 (1):\penalty0
  258--265, 1996{\natexlab{b}}.

\bibitem[Pardoux and Peng(1992)]{PaPe:92}
E.~Pardoux and S.~Peng.
\newblock Backward stochastic differential equations and quasilinear parabolic
  partial differential equations.
\newblock In \emph{Stochastic Partial Differential Equations and Their
  Applications}, pages 200--217. Springer, 1992.

\bibitem[Pardoux and Tang(1999)]{PaTa:99}
E.~Pardoux and S.~Tang.
\newblock Forward-backward stochastic differential equations and quasilinear
  parabolic {PDEs}.
\newblock \emph{Probability Theory and Related Fields}, 114\penalty0
  (2):\penalty0 123--150, 1999.

\bibitem[Prasad and Sethi(2004)]{PrSe:04}
A.~Prasad and S.~P. Sethi.
\newblock Competitive advertising under uncertainty: A stochastic differential
  game approach.
\newblock \emph{Journal of Optimization Theory and Applications}, 123\penalty0
  (1):\penalty0 163--185, 2004.

\bibitem[Robinson(1951)]{Ro:51}
J.~Robinson.
\newblock An iterative method of solving a game.
\newblock \emph{Annals of Mathematics}, pages 296--301, 1951.

\bibitem[Sannai et~al.(2019)Sannai, Takai, and Cordonnier]{Sannai2019universal}
A.~Sannai, Y.~Takai, and M.~Cordonnier.
\newblock Universal approximations of permutation invariant/equivariant
  functions by deep neural networks.
\newblock \emph{arXiv preprint arXiv:1903.01939}, 2019.

\bibitem[Seale and Burnett(2006)]{SeBu:06}
D.~A. Seale and J.~E. Burnett.
\newblock Solving large games with simulated fictitious play.
\newblock \emph{International Game Theory Review}, 8\penalty0 (03):\penalty0
  437--467, 2006.

\bibitem[Shapley(1964)]{Sh:64}
L.~S. Shapley.
\newblock Some topics in two-person games.
\newblock \emph{Advances in Game Theory}, 52:\penalty0 1--29, 1964.

\bibitem[Van~Long(2011)]{Va:11}
N.~Van~Long.
\newblock Dynamic games in the economics of natural resources: a survey.
\newblock \emph{Dynamic Games and Applications}, 1\penalty0 (1):\penalty0
  115--148, 2011.

\bibitem[Zaheer et~al.(2017)Zaheer, Kottur, Ravanbakhsh, Poczos, Salakhutdinov,
  and Smola]{ZaheerKotturRavanbakhshEtAl2017}
M.~Zaheer, S.~Kottur, S.~Ravanbakhsh, B.~Poczos, R.~R. Salakhutdinov, and A.~J.
  Smola.
\newblock Deep sets.
\newblock In \emph{Advances in Neural Information Processing Systems}, pages
  3391--3401, 2017.

\bibitem[Zhang et~al.(2018{\natexlab{a}})Zhang, Yang, Liu, Zhang, and
  Basar]{zhang2018fully}
K.~Zhang, Z.~Yang, H.~Liu, T.~Zhang, and T.~Basar.
\newblock Fully decentralized multi-agent reinforcement learning with networked
  agents.
\newblock In \emph{Proceedings of the 35th International Conference on Machine
  Learning}, pages 5872--5881, 2018{\natexlab{a}}.

\bibitem[Zhang et~al.(2018{\natexlab{b}})Zhang, Han, Wang, Car, and
  E]{ZhHaE:18}
L.~Zhang, J.~Han, H.~Wang, R.~Car, and W.~E.
\newblock Deep potential molecular dynamics: a scalable model with the accuracy
  of quantum mechanics.
\newblock \emph{Physical Review Letters}, 120\penalty0 (14):\penalty0 143001,
  2018{\natexlab{b}}.

\bibitem[Zhang et~al.(2018{\natexlab{c}})Zhang, Han, Wang, Saidi, Car, and
  E]{Zhang2018end}
L.~Zhang, J.~Han, H.~Wang, W.~Saidi, R.~Car, and W.~E.
\newblock End-to-end symmetry preserving inter-atomic potential energy model
  for finite and extended systems.
\newblock In \emph{Advances in Neural Information Processing Systems}, pages
  4436--4446, 2018{\natexlab{c}}.

\end{thebibliography}

\newpage
\appendix
\section{Technical Details to Numerical Examples}\label{appendix_PDE}
\subsection{The Analytic Solution in Section~\ref{sec_example1}}
\label{appendix_example1}
The results in this section is firstly derived in \cite[Section~3]{CaFoSu:15}, and we repeat them here for completeness. Assume the following ansatz for the HJB system \eqref{eq_HJBexample1}:
\begin{equation}
V^i(t,\bm x) = \frac{\eta(t)}{2}(\bar x - x^i)^2 + \mu(t), \quad i \in \mc{I},
\end{equation}
where $\eta(t), \mu(t)$ are two scalar functions to be determined.
Under this ansatz the optimal feedback control becomes
\begin{equation}\label{eq:example1_alpha}
\alpha^{i,\ast}(t, \bm x) = \left[q + \eta(t)(1-\oon)\right](\bar x - x^i).
\end{equation}
Plugging the ansatz into \eqref{eq_HJBexample1} and collecting the coefficients of the squared and constant terms, we find that $\eta(t)$ solves a Riccati equation
\begin{align}
\dot \eta(t) &= 2(a+q)\eta(t) + (1-1/N^2) \eta^2(t) - (\eps - q^2), \quad \eta(T) = c,
\end{align}
and $\mu(t)$ depends on $\eta(t)$ through:
\begin{align}
\dot \mu(t) = -\half \sigma^2(1-\rho^2)(1-1/N)\eta(t), \quad \mu(T) = 0.
\end{align}
The solution to the Riccati equation is 
\begin{align}
\eta(t) &= \frac{-(\eps - q^2)(e^{(\delta^+ - \delta^-)(T-t)}-1) - c (\delta^+e^{(\delta^+ - \delta^-)(T-t)} - \delta^-)}{(\delta^-e^{(\delta^+ - \delta^-)(T-t)}-\delta^+) - c(1-1/N^2)(e^{(\delta^+ - \delta^-)(T-t)}-1)},
\end{align}
where $
\delta^{\pm} = -(a+q) \pm \sqrt R$, and $R = (a+q)^2 + (1-1/N^2)(\eps - q ^2).
$

\subsection{Technical Details in Section~\ref{sec_example2}}\label{appendix_example2}
In this risk-sensitive dynamic game, the Nash equilibrium is related to the following system:
\begin{align}
\label{eq:example2_appendix}
\begin{dcases}
\partial_t V^i + \inf_{\alpha^i}\left\{\sum_{j=1}^N [a(\bar x - x^j) + \alpha^j] \partial_{x^j}V^i + V^i\theta_i\left(\frac{(\alpha^i)^2}{2} - q\alpha^i(\bar x - x^i) + \frac{\eps}{2}(\bar x - x^i)^2 \right)\right\} \\
\quad\quad\quad + \half \text{Tr}(\Sigma\transpose \text{Hess}_{\bx} V^i \Sigma) = 0,\\
V^i(T,\bm x) =  \theta_i \exp\left\{\theta_i \frac{c}{2}(\bar x - x^i)^2\right\},
\end{dcases}
\end{align} 
and the candidate of player $i$'s optimal strategy is given by
\begin{equation}\label{eq:example2_alpha}
\alpha^{i}(t, \bm x) = q(\bar x - x^i)- \frac{\partial_{x^i} V^i(t, \bx)}{\theta_i V^i(t, \bx)}.
\end{equation}

Observing the exponential quadratic form in the terminal condition and in the cost functional, we assume the following ansatz
\begin{equation}
V^i(t,\bm x) = \theta_i \exp\left\{ \theta_i \left(\half \bm x\transpose P^i(t) \bm x + p^i(t)\right)\right\}, \quad i\in \mc{I},
\end{equation}
where $P^i(t) \in \RR^{n\times n}, p^{i}(t)\in \RR$ are two (matrix) functions to be determined.
Under this ansatz the optimal feedback control becomes
\begin{equation}
\alpha^{i,\ast}(t, \bm x) = q(\bar x - x^i) - P^i(t)\bx.
\end{equation}
Plugging the ansatz into \eqref{eq:example2_appendix} and collecting the coefficients of the squared and constant terms, 
we find that $P^i(t)$ and $p^i(t)$ satisfy
\begin{equation}
\dot P^i + 2(a+q) A P^i - 2 \sum_{j=1}^N (P^j)\transpose \Delta_{j,j} P^i + (P^i)\transpose (\Delta_{i,i} + \theta_i \Sigma\Sigma\transpose) P^i + (\eps - q^2)e_i e_i\transpose = 0, \quad P^i(T) = ce_i e_i^T,
\end{equation}
and 
\begin{equation}
\dot p^i(t) + \half \text{Tr}(\Sigma\transpose P^i(t)\Sigma) = 0, \quad p^i(T) = 0.
\end{equation}
Here we have used the notation
\begin{align}
&A = (A_{i,j}) = (\oon - \dij{i,j}), \quad e_i = \Big[\oon, \ldots, \oon, \place{i^{th}}{\oon-1}, \oon, \ldots, \oon\Big]\transpose,\\
& B_i = [0,\ldots, 0, \place{i^{th}}{1}, 0, \ldots, 0]\transpose,  \quad  \Delta_{i,i} = B_i B_i\transpose,
\end{align}
and $\Sigma$ is defined as in \eqref{eq:example1_driftdiff}. The solutions of $P^i(t)$ and $p^i(t)$ then can be obtained through high-precision numerical integration. Note that these matrix Riccati equations for $\{P^i(t), i \in \mc{I}\}$ are coupled, while $p^i(t)$ depends solely on $P^i(t)$.

To obtain the HJB equation for the individual decision problem decoupled by fictitious play, we plugin \eqref{eq:example2_alpha} into \eqref{eq:example2_appendix} and deduce
\begin{align}
\partial_t V^{i} + \half \text{Tr}(\Sigma\transpose \text{Hess}_{\bx} V^i \Sigma) + a(\bar x - x^i)\partial_{x^i}V^{i} + \sum_{j\neq i} [a(\bar x - x^j) + \alpha^{j}(t,\bm x)] \partial_{x^j}V^{i} \\
\quad\quad\quad + \theta_i V^i \left(\frac\eps2 (\bar x - x^i)^2 - \half \left(q(\bar x - x^i) - \frac{\partial_{x^i}V^i}{\theta_i V^i}\right)^2\right) = 0.
\end{align} 
This is indeed of the form \eqref{eq:FP_PDE_explicit}, with $\mu^i(t, \bm x; \balpha^{-i})$ following the definition in Section~\ref{sec_example1}, and $h^i$ given by:
\begin{equation}
    h^i(t,\bm x, y, \bm z; \balpha^{-i}) = \theta_i y\left(\frac{\eps}{2}(\bar x - x^i)^2 - \half (q(\bar x - x^i) - \frac{z^i}{\theta_i y \sigma \sqrt{1-\rho^2}} )^2\right).
\end{equation}

\subsection{Technical Details in Section~\ref{sec_example3}}\label{appendix_example3}
In this general linear-exponential-quadratic game, the coupled HJB system can be written as
\begin{align}\label{eq:example3}
\begin{dcases}
V_t^i + \half \text{Tr}(\Sigma\transpose \text{Hess}_{\bx}V^i\Sigma) + \inf_{\alpha^i \in \RR^{n_i}} \left\{\nabla_{\bx} V^i \cdot (A\bm x + \sum_{j=1}^N B_j \alpha^j) +  \frac12 \theta_i V^i\left(\bx\transpose Q_i \bx + (\alpha^i)\transpose R_i \alpha^i\right)\right\}  = 0, \\
V^i(T,\bm x) = \theta_i \exp\left\{\frac{\theta_i}{2}\bm x\transpose M_i \bm x\right\}, \quad i \in \mc{I}. 
\end{dcases}
\end{align} 
The optimal control satisfies
\begin{equation}
\alpha^{i}(t, \bm x) = - \frac{R_i\inv B_i\transpose \nabla_{\bx} V^i(t, \bx)}{\theta_i V_i(t, \bx)}.
\label{eq:example3_alpha}
\end{equation}
In order to solve the system \eqref{eq:example3} directly for the Nash equilibrium, like in the second example, we assume the following ansatz with the exponential form
\begin{equation}
V^i(t,\bm x) = \theta_i \exp\left\{ \theta_i \left(\half \bm x\transpose P^i(t) \bm x + p^i(t)\right)\right\}, \quad i\in \mc{I},
\end{equation}
where $P^i(t) \in \RR^{n\times n}, p^{i}(t)\in \RR$ are two functions of $t$ to be determined.
Using this ansatz, the optimal control becomes
\begin{equation}
\alpha^{i,\ast}(t,\bm x) = - R_i\inv B_i\transpose P^i(t)\bm x.
\end{equation}
Plugging the ansatz into \eqref{eq:example3}, we deduce that $P^i(t)$ solves a matrix Riccati equation
\begin{align}\label{eq_P_matrix}
&\dot P^i + (P^i)\transpose A + A\transpose P^i + Q_i + (P^i)\transpose(B_i R_i\inv B_i\transpose + \theta_i \Sigma\Sigma\transpose)P^i \\
&\quad\quad\quad\quad\quad - (P^i)\transpose \sum_{j=1}^N B_j R_j\inv B_j\transpose P^j - ( \sum_{j=1}^N B_j R_j\inv B_j\transpose P^j)\transpose P^i = 0, \quad P^i(T) = M_i,
\end{align}
and $p^i(t)$ solves the ODE
\begin{equation}
\dot p^i(t) + \half\text{Tr}(\Sigma\transpose P^i(t)\Sigma) = 0, \quad p^i(T) = 0.
\end{equation}
We remark that these Riccati equations are coupled as well, and can be solved by high-precision numerical integration.

To obtain the HJB equation for the individual decision problem decoupled by fictitious play,
we plugin\eqref{eq:example3_alpha} into \eqref{eq:example3} and deduce the simplified PDE
\begin{multline}
    V_t^i + \half \text{Tr}(\Sigma\transpose \text{Hess}_{\bx}V^i\Sigma) + \nabla_{\bx} V^i \cdot (A\bm x + \sum_{j \neq i} B_j \alpha^j) + \frac12 \theta_i V^i \bx\transpose Q_i \bx \\
    - \frac12 \frac{(\nabla_{\bx} V^i)\transpose B_i R_i\inv B_i\transpose \nabla_{\bx} V^i}{\theta_i V^i} = 0.
\end{multline}
This PDE has the same form of \eqref{eq:FP_PDE_explicit}, with
\begin{align}
    \mu^i(t,\bm x; \balpha^{-i}) &=  A \bm x + \sum_{j \neq i} B_j \alpha^{j}(t, \bm x),\\
    h^i(t, \bm x, y, \bm z; \balpha^{-i}) &= \frac{1}{2}\theta_i y \bm x\transpose Q_i \bm x - \frac{1}{2} \frac {\bm z\transpose \Sigma^{-1}B_i R_i\inv B_i\transpose(\Sigma^{-1})\transpose \bm z}{\theta_i y}.
\end{align}

\section{Hyperparameters and Details of Deep Fictitious Play}\label{appendix_parameter}
This appendix reports the hyperparameters used in the numerical examples of Section~\ref{sec_numerics}. We fix $N_{\text{SGD\_per\_stage}}=100$ and set $N_{\text{sample}}=N_{\text{SGD\_per\_stage}} N_{\text{batch}}$. This means that the samples are never reused in the learning. In all the numerical examples, the value functions $V^i$ are approximated by feedforward neural networks $\text{Net}(t, \bx)$ with 3 hidden layers whose widths are all the same. 
We use {\it{tanh}} as the activation function and adopt batch normalization (\cite{batch}) right after each affine transformation and before activation.
Each component of $\bx_0$, the initial state $\bX_0^i$, is sampled independently from the uniform distribution on $[-\delta_0, \delta_0]$. $\delta_0$ is chosen such that in the following process driven by the optimal policy $\balpha^*$
\begin{equation*}
\ud \bm X_t^{\bm \alpha^*} = b(t, \bm X_t^{\bm \alpha^*}, \bm \alpha^*(t, \bm X_t^{\bm \alpha^*})) \ud t + \Sigma(t, \bm X_t^{\bm \alpha^*}) \ud \bm W_t, \quad \bm X_0 = \bx_0,
\end{equation*}
the standard deviation of $\{\bX_t\}_{t=0}^T$ is approximately $\delta_0$. In other words, $\delta_0$ is determined as a fixed-point. 
The rationale for such a procedure is to make sure the data generated for the learning is representative enough in the whole state space. We also test other reasonable choices of $\bx_0$ and find the final accuracy is insensitive to it. For the inter-bank game in Section~\ref{sec_example4} with $N=50$ and the superlinear dynamics~\eqref{eq:ex4_dynamics}, the optimal policy is unavailable for determining $\delta_0$. We instead use $\delta_0$ defined for the game with $N=50$ and the linear dynamics~\eqref{eq:ex1_dynamics}.

Table~\ref{tab_parameter} reports the values of some other hyperparameters used in the numerical examples. 
\begin{table}[h]
\caption{Hyperparameters and runtime for the numerical examples presented in Section 4.}\label{tab_parameter}
\begin{center}
~\newline 
\begin{tabular}{@{}c|cccc@{}}
\toprule
Parameters / Problem & Section 4.1     & Section 4.2 & Section 4.3 & Section 4.4\\ \midrule
$N_T$       &  40 &  40 &  30  & 40 \\
width of hidden layers &  40 &  40 &  40  & 60 \\
$M$ (\# of total stages)       &  80 &  100 &  200  & 400\\
$N_\text{batch}$ &  256 &  512 &  512 & 256 \\
learning rate     & 5e-4 &  5e-4  &  (5e-3, 5e-4)$^\ast$ & 5e-4\\ 
runtime (hours) $^\dagger$ & 7 & 13 & 20.5 & 35\\
\bottomrule
\end{tabular}
\end{center}
\small{$^\ast$The learning rate is piecewise constant in the example of Section 4.3. It equals 5e-3 for the first half of SGD updates and 5e-4 for the second half.\newline
$^\dagger$ The numerical experiments were conducted on an NVIDIA Tesla P100 GPU. The runtime is subject to further reduction with a multi-GPU system.
}
\end{table}

In Section~\ref{sec_example4} we solve the large scale multi-agent game with $N=50$. Due to the intrinsic symmetry of the game, we assume all the players essentially share the same strategies during each stage of the fictitious play. Without loss of generality, we only keep track of player 1's strategy and denote its neural network approximation by
\begin{equation}
\label{single_alpha}
    \alpha^{1,m}(t, x_1, x_2, \dots, x_N) = \text{Net}(t, x_1, x_2, \dots, x_N; \bm{w}),
\end{equation}
where $\text{Net}$ is a neural network and $\bm w$ denotes all the trainable parameters. Then player $i$'s strategy $(i\neq 1)$ is defined by swapping the components $x_1$ and $x_i$ in the arguments, i.e.,
\begin{equation}
\label{symmetry_alpha}
    \alpha^{i,m}(t, x_1, x_2, \dots, x_N) = \text{Net}(t, x_i, \dots, x_{i-1}, x_1, x_{i+1}, \dots, x_N; \bm{w}).
\end{equation}
This completes the simplification of Algorithm~\ref{def_algorithm1}, whose pseudo-code is summarized in Algorithm~\ref{def_algorithm2} below. 

Note that there are other ways besides~\eqref{symmetry_alpha} to define the others' strategies based on~\eqref{single_alpha}, for instance,
\begin{equation}
\label{symmetry_alpha2}
    \alpha^{i,m}(t, x_1, x_2, \dots, x_N) = \text{Net}(t, x_i, \dots, x_{i+1}, x_1, x_{i-1}, \dots, x_N; \bm{w}).
\end{equation}
Due to the total symmetry of this game, the optimal strategy $\alpha^{1,*}$ should be permutation invariant with respect to the arguments $x_2, \dots, x_N$, which implies that~\eqref{symmetry_alpha} and ~\eqref{symmetry_alpha2} are the same if $\text{Net}(\cdot; \bm{w})$ represents $\alpha^{1,*}$ exactly. In the current algorithm, the plain feedforward neural network does not guarantee this property, but we find different choices (like ~\eqref{symmetry_alpha} and ~\eqref{symmetry_alpha2}) do not make a substantial difference in final accuracy. On the hand, designing specific neural network architectures to guarantee permutation invariance/equivariance exactly has been discussed in depth in the recent literature (\cite{ZaheerKotturRavanbakhshEtAl2017,Zhang2018end,Sannai2019universal}). It will be of interest to test in future work whether such networks can improve the performance in learning the Nash equilibrium of totally symmetric games.

\begin{algorithm}[H]
\caption{Deep Fictitious Play for Finding Markovian Nash Equilibrium of Symmetric Game \label{def_algorithm2}}
    \begin{algorithmic}[1]
	\REQUIRE $N$ = \# of players, $N_T$ = \# of subintervals on $[0,T]$, $M$ = \# of total stages in fictitious play, $N_{\text{sample}}$ = \# of sample paths generated for each player at each stage of fictitious play, $N_{\text{SGD\_per\_stage}}$ = \# of SGD steps for each player at each stage, $N_{\text{batch}}$ = batch size per SGD update, $\alpha^{1, 0}\colon$ the initial smooth policy for player 1
	    \STATE Initialize one deep neural network to represent $V^{1,0}$ for player 1
	    \STATE Define $\alpha^{2, 0}, \dots, \alpha^{N, 0}$ according to \eqref{symmetry_alpha} and collect $\balpha^0 \gets (\alpha^{1,0}, \ldots, \alpha^{N,0})$
		\FOR{$m \gets 1$ to $M$}
		\STATE  Generate $N_\text{sample}$ sample paths $\{\bX_{k}^{1,\pi}\}_{k=0}^{N_T}$ according to \eqref{eq:disc_X_path} and the realized optimal policies $\balpha^{-1, m-1}(t_k, \bm X_k^{1,\pi})$
		\FOR{$\ell \gets 1$ to $N_{\text{SGD}\_\text{per}\_\text{stage}}$}
		    \STATE Update the parameters of the neural network one step with $N_{\text{batch}}$ paths using the SGD algorithm (or its variant), based on the loss function \eqref{eq:disc_objective}
		  \ENDFOR
		  \STATE Obtain the approximate optimal policy  $\alpha^{1,m}$  according to \eqref{def_alphaast}
		  \STATE Define $\alpha^{2, m}, \dots, \alpha^{N, m}$ according to \eqref{symmetry_alpha} and collect the optimal policies at stage $m$:\\
		  $\bm{\alpha}^m \gets (\alpha^{1,m}, \ldots, \alpha^{N,m})$
		  \ENDFOR
		\RETURN The optimal policy $\balpha^{M}$
	\end{algorithmic}
\end{algorithm}

\section{Robustness of the Deep Fictitious Play}\label{appendix_sensitivity}

\begin{figure}[!ht]
        \centering
        \includegraphics[width=0.95\textwidth]{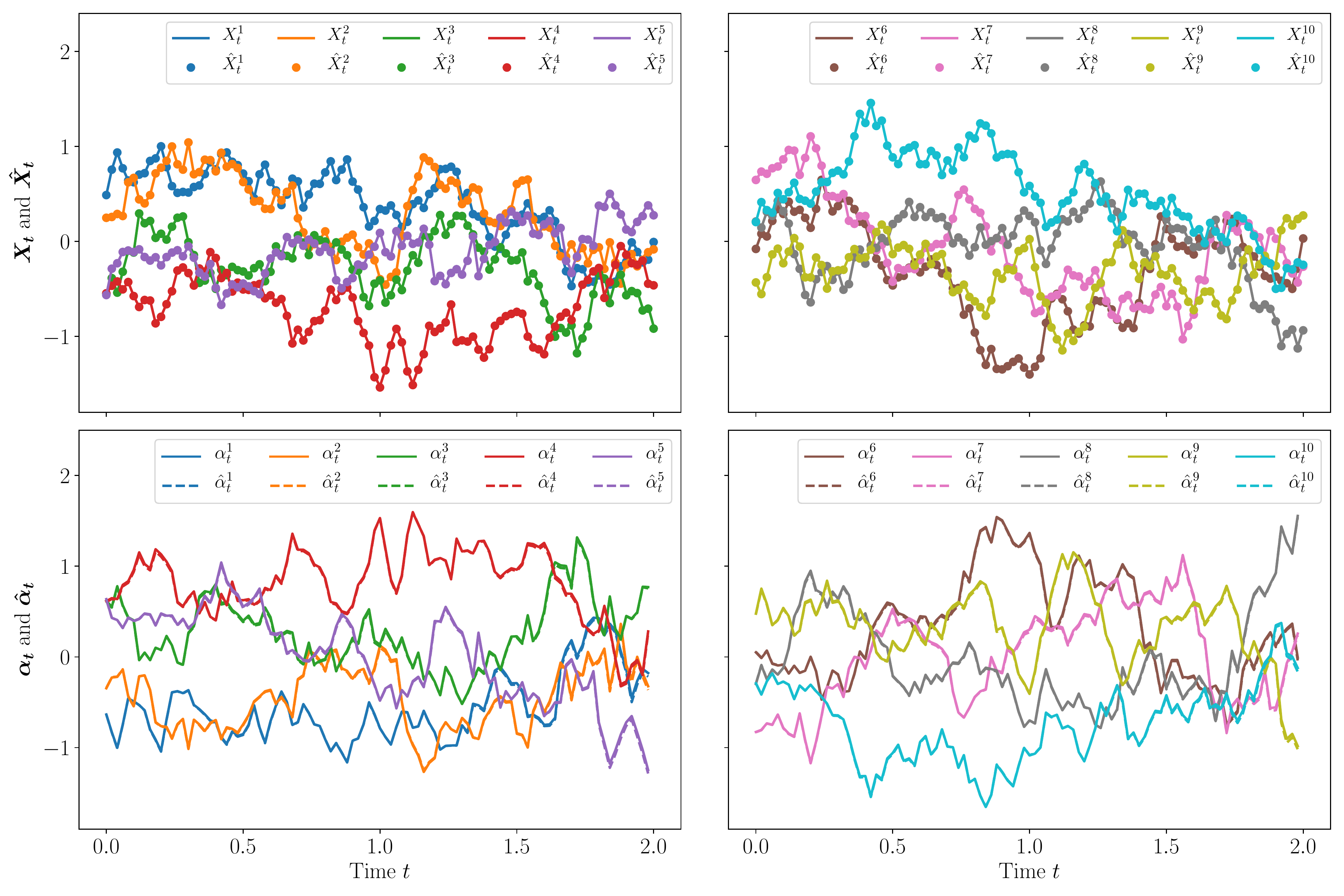}
        \caption{A sample path for each player of the inter-bank game in Section 4.1 with $N=10$ and the parameters in Eq.~\eqref{def_parameters} except $T=2$. Top: the optimal state process $X_t^i$ (solid lines) and its approximation $\hat{X}_t^i$ (circles) provided by the optimized neural networks, under the same realized path of Brownian motion. Bottom: comparisons of the strategies $\alpha_t^i$ and $\hat{\alpha}_t^i$ (dashed lines).}
        \label{fig:LQMFG_path_t2}
    \end{figure}

\begin{figure}[!ht]
        \centering
        \includegraphics[width=0.95\textwidth]{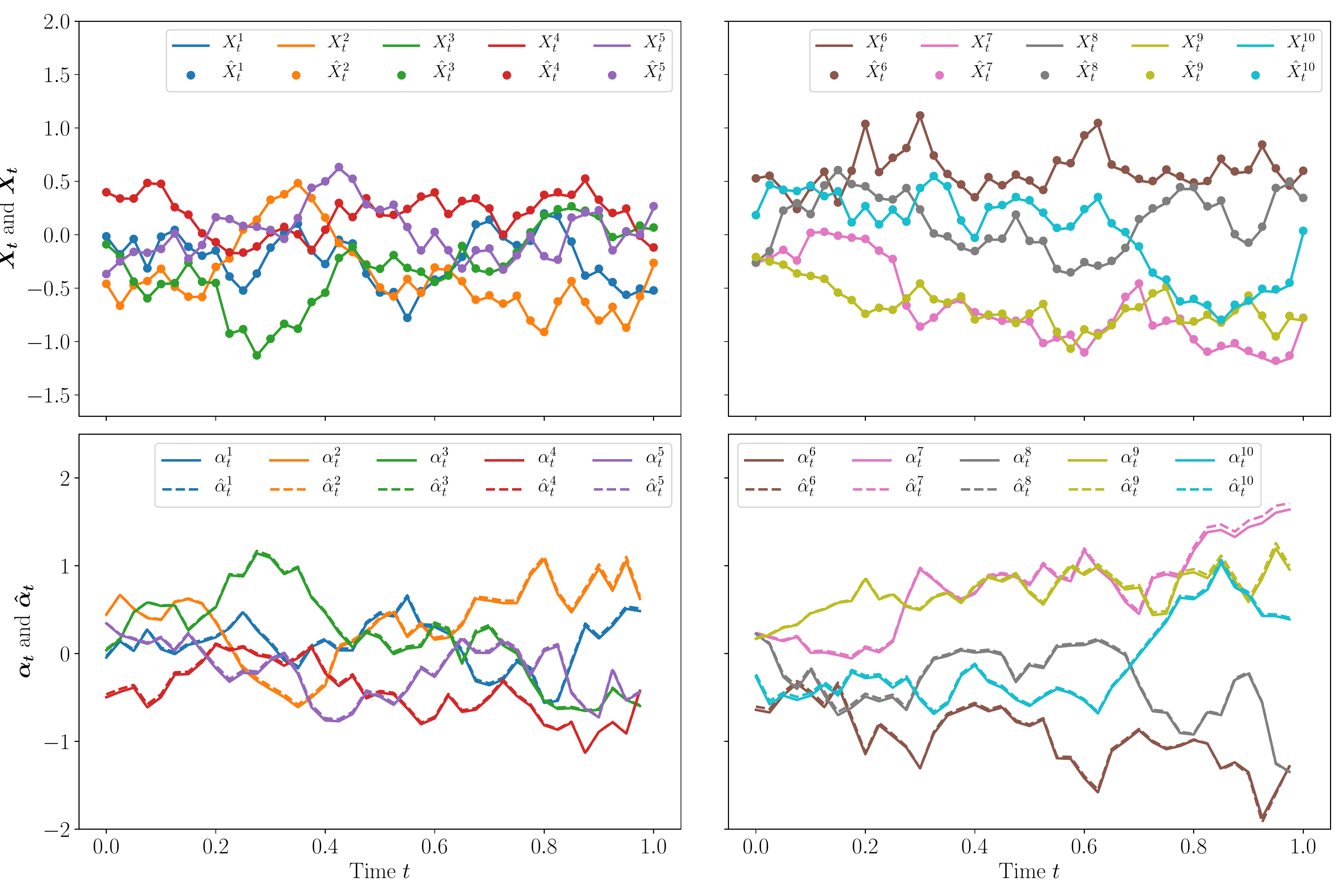}
        \caption{A sample path for each player of the inter-bank game in Section 4.1 with $N=10$  and the parameters in Eq.~\eqref{def_parameters} except $a=0.3$. Top: the optimal state process $X_t^i$ (solid lines) and its approximation $\hat{X}_t^i$ (circles) provided by the optimized neural networks, under the same realized path of Brownian motion. Bottom: comparisons of the strategies $\alpha_t^i$ and $\hat{\alpha}_t^i$ (dashed lines).}
        \label{fig:LQMFG_path_a0p3}
    \end{figure}

Regarding robustness of the deep fictitious play, we did two more tests on the inter-bank game in Section~\ref{sec_example1}. Compared to the original parameter (see Eq. \eqref{def_parameters}), we solve two modified cases: (1) increasing the game length $T=2$ (with $N$ increased to 100) instead of $T=1$, and (2) increasing the coupling $a=0.3$ instead of $a=0.1$. Figures~\ref{fig:LQMFG_path_t2} and \ref{fig:LQMFG_path_a0p3} present the sample paths for each player of the optimal state process $X_t^i$ and the optimal control $\alpha_t^i$ \emph{vs.} their approximations $\hat{X}_t^i, \hat{\alpha}_t^i$ provided by the optimized neural networks. In these two cases, the final RSE of $V$ and $\nabla V$ are (1) 4.1\% and 0.2\%, (2) 4.9\% and 0.1\%. We have not identified bifurcation in these experiments, which exists in algorithms of solving mean field games (cf. \cite{AnCh:18}). In other words, the performance of deep fictitious play should be robust at least to a certain range of parameters.

\end{document}